\documentclass[12pt]{amsart}
\usepackage{amssymb}
\usepackage[cmtip,matrix,arrow]{xy}

\newtheorem{theorem}{Theorem}[section]

\newtheorem{proposition}[theorem]{Proposition}
\newtheorem{definition}[theorem]{Definition}

\newtheorem{lemma}[theorem]{Lemma}
\theoremstyle{definition}    
\theoremstyle{remark}

\newtheorem{remark}[theorem]{Remark}

\newtheorem{example}[theorem]{Example}
\newtheorem{examples}[theorem]{Examples}

\newcommand\A{\mathcal{A}}
\newcommand\M{\mathcal{M}}

\newcommand{\W}{\mathcal{W}}

\renewcommand{\L}{\mathcal{L}}
\renewcommand{\O}{\mathcal{O}}

\newcommand{\Co}{\mathcal{C}}
\newcommand{\ca}{\mathcal}

\newcommand{\U}{\on{U}}

\newcommand{\R}{\mathbb{R}}
\newcommand{\C}{\mathbb{C}}
\newcommand{\SU}{\on{SU}}

\newcommand{\Z}{\mathbb{Z}}

\newcommand\lie[1]{\mathfrak{#1}}

\newcommand{\h}{\lie{h}}
\newcommand{\g}{\lie{g}}
\newcommand{\m}{\lie{m}}

\renewcommand{\t}{\lie{t}}

\newcommand{\on}{\operatorname}

\newcommand{\Ad}{ \on{Ad} }
\newcommand{\ad}{\on{ad}}

\renewcommand{\ker}{ \on{ker}}

\newcommand{\Spin}{ \on{Spin}}
\newcommand{\Pin}{ \on{Pin}}
\newcommand{\SO}{ \on{SO}}
\newcommand{\Mult}{  \on{Mult}}
\newcommand{\Vol}{  \on{Vol}}
\newcommand{\vol}{  \on{vol}}
\newcommand{\diag}{  \on{diag}}
\renewcommand{\dh}{\m}

\renewcommand{\o}{{\on{o}}}

\newcommand\qu{/\kern-.7ex/} 
\newcommand{\fus}{\circledast} 

\newcommand{\hra}{\hookrightarrow}

\renewcommand{\d}{{\mbox{d}}}
\newcommand{\ol}{\overline}

\newcommand\Phinv{\Phi^{-1}}

\newcommand\lam{\lambda}
\newcommand\Lam{\Lambda}
\newcommand\Sig{\Sigma}
\newcommand\sig{\sigma}
\newcommand\eps{\epsilon}
\newcommand\Om{\Omega}
\newcommand\om{\omega}
\newcommand{\Del}{\Delta}

\newcommand{\f}{\frac}

\newcommand{\p}{\partial}
\renewcommand{\l}{\langle}
\renewcommand{\r}{\rangle}
\newcommand\hh{{\f{1}{2}}}
\newcommand{\ti}{\tilde}
\newcommand{\eeq}{\end{eqnarray*}}
\newcommand{\beq}{\begin{eqnarray*}}

\renewcommand{\H}{\ca{H}}
\newcommand{\La}{\Lambda}

\newcommand{\Cl}{{\on{Cl}}}

\newcommand{\wh}{\widehat}
\newcommand{\wt}{\widetilde}
\newcommand{\mf}{\mathfrak}
\begin{document}
\sloppy

\title[Duistermaat-Heckman measure]
{Duistermaat-Heckman measures and \\
moduli spaces of flat bundles over surfaces}

\author{A. Alekseev}
 \address{Institute for Theoretical Physics \\ Uppsala University \\
 Box 803 \\ \mbox{S-75108} Uppsala \\ Sweden}
 \email{alekseev@teorfys.uu.se}

\author{E. Meinrenken}
 \address{University of Toronto, Department of Mathematics,
 100 St George Street, Toronto, Ontario M5S3G3, Canada }
 \email{mein@math.toronto.edu}

\author{C. Woodward}
 \address{Mathematics-Hill Center, Rutgers University,
110 Frelinghuysen Road, Piscataway NJ 08854-8019, USA}
\email{ctw@math.rutgers.edu}

\date{February 1999. Revised, April 2001}
\begin{abstract}
We introduce Liouville measures and Duistermaat-Heckman measures for
Hamiltonian group actions with group valued moment maps. The theory is
illustrated by applications to moduli spaces of flat bundles on
surfaces.
\end{abstract}
\maketitle
%

\section{Introduction}

One of the fundamental invariants of a Hamiltonian $G$-manifold $M$ in
symplectic geometry is the Duistermaat-Heckman (DH) measure on the
dual of the Lie algebra $\g^*$, defined as the push-forward of the
Liouville measure under the moment map. The DH-measure encodes volumes
of reduced spaces, and by the Duistermaat-Heckman theorem \cite{du:on}
its derivatives describe Chern numbers for the corresponding level
set.

The purpose of this paper is to develop a theory of Liouville measures
and DH measures for Hamiltonian $G$-manifolds $(M,\om,\Phi)$ with
group valued moment maps $\Phi:\,M\to G$ as introduced in
\cite{al:mom}.  The basic examples of such spaces are conjugacy
classes $\Co\subset G$, with moment map the inclusion into $G$. 
Since we allow $G$ to be disconnected, these include 
all compact symmetric spaces (up to finite covers).  Another example is
the space $G^{2h}$, with moment map the product of Lie group
commutators, $(a_1,b_1,\ldots,a_h,b_h)\mapsto \prod_j [a_j,b_j]$.
There is a notion of reduction for group valued moment maps, and the
reduced spaces $M\qu G=\Phinv(e)/G$ are symplectic. For instance,
$G^{2h}\qu G$ is the moduli space of flat $G$-bundles over a surface
of genus $h$, with its natural symplectic structure
\cite{at:ge,at:mom}.

The 2-form $\om$ for a space with group valued moment map is usually 
degenerate. The kernel of $\om$ is controlled by the 
{\em minimal degeneracy} condition, involving the moment map $\Phi$. 
If the group $G$ is 1-connected, we obtain a volume form 
as the top degree part of 
$$\Gamma=e^\om\,\Phi^*\ca{T},$$ 
where $\ca{T}$ is a differential form on $G$ constructed from the
homomorphism $G\to \on{Spin}(\g)$. In \cite{al:cl}, the formula for 
the volume form is heuristically obtained as a Feynman integral 
over the path fibration $P_0G\to G$. For general compact groups $G$,
group valued Hamiltonian $G$-manifolds are not necessarily orientable:
an example is the conjugacy class in $\SO(3)$ consisting of rotations
by $\pi$, which is isomorphic to the real projective plane. In such
cases, we can still define $\Gamma$ as a form with values in the
orientation bundle, and obtain a nowhere vanishing density on $M$. We
call this density the Liouville measure, and its push-forward to $G$
the DH measure of $(M,\om,\Phi)$.  Just as in the $\g^*$-valued
theory, the DH-measures describe volumes of reduced spaces. One of our
main results says that the Liouville measure for a product of two
spaces with group valued moment map is equal to the direct product of
the Liouville measures, and consequently the DH-measure is the
convolution on $G$ of the DH-measures of its factors.

For the example $G^{2h}$, we show that the Liouville measure coincides
with Haar measure.  The push-forward of Haar measure on $G^{2h}$ under
the map $(a_1,b_1,\ldots,a_h,b_h)\mapsto \prod_j [a_j,b_j]$ gives
Witten's Fourier series \cite[Equation (2.73)]{wi:qg}.  Witten,
Bismut-Labourie \cite{bi:sy} and Liu \cite{li:he} use a Reidemeister
torsion calculation to identify this expression with the symplectic
volume of moduli spaces of flat $G$-bundles over surfaces.  Our
construction provides a purely symplectic approach to the Witten
volume formulas, and gives an extension to disconnected compact Lie
groups.  

\vskip.1in

\noindent{\bf Acknowledgment:} 
We thank N. Berline and B. Kostant for discussion, and M. Vergne for a
number of helpful comments on this paper. 
\vskip.3in

\begin{center}
{\bf Notation.}
\end{center}
\vskip.1in

Throughout this paper $G$ denotes a compact  Lie group, 
and $G^\o$ its identity component. 
The left/right invariant 
Maurer-Cartan forms on $G$ are denoted $\theta^L,\theta^R\in\Om^1(G,\g)$. 
We fix an invariant inner 
product on the Lie algebra $\g$ of $G$. The corresponding Riemannian 
measure on $G$ will be denoted $\d\vol_G$, and the Riemannian 
volume $\vol G$. 

A $G$-manifold is a manifold $M$ together with a group homomorphism
$\phi:\,G\to \on{Diff}(M)$ such that the action map $G\times M\to
M,\,(g,m)\mapsto g.m:=\phi(g)(m)$ is smooth.  We denote by $G_m$ the
stabilizer group of $m\in M$ and by $\g_m$ its Lie algebra.  Given
$\xi\in\g$, we denote by $\xi_M$ the corresponding generating vector
field: $\xi_M(m)=\f{d}{d\,t}|_{t=0}\exp(-t\xi).m$. Let 
$\Om_G(M)$ be the Cartan complex of equivariant differential forms
on $M$. Thus $\Om_G(M)$
is the space of equivariant maps $\alpha:\,\g\to \Om(M)$, with
equivariant differential $(\d_G\alpha)(\xi)=
(\d+\iota(\xi_M))\alpha(\xi)$.

\section{Group valued Hamiltonian $G$-manifolds}
In this Section we recall the definition of a space with group valued 
moment map, and describe the main examples. 
We refer to \cite{al:mom} for proofs and further details. 

\subsection{Definition of a group valued Hamiltonian $G$-manifold}
Let $G$ be a compact Lie group, acting on itself by conjugation. 
The canonical closed 3-form 
\begin{equation}
\eta=\f{1}{12}\theta^L\cdot[\theta^L,\theta^L]
\end{equation}
on $G$ extends to an equivariantly closed equivariant 3-form, 
\begin{equation} 
\eta_G(\xi)=\eta+\f{1}{2}(\theta^L+\theta^R)\cdot\xi.
\end{equation}
\begin{definition}\label{def:qham}
A group valued Hamiltonian $G$-manifold (called quasi-Hamiltonian in
\cite{al:mom}) is a triple $(M,\om,\Phi)$
consisting of a $G$-manifold $M$, an invariant 2-form $\om$, and an
equivariant map $\Phi\in C^\infty(M,G)^{G}$ such that 
\begin{equation}\label{eq:cond1}
\d_G\om=\Phi^*\eta_G
\end{equation}
and such that for all $m\in M$, $\ker\om_m=\{\xi_M(m)|\,
(\on{Ad}_{\Phi(m)}+1)\xi=0\}.$
\end{definition}
Condition \eqref{eq:cond1} splits into 
$\d\om=\Phi^*\eta$ and the {\em moment map condition}
\begin{equation}\label{eq:cond3}
\iota(\xi_M)\om=\hh \Phi^*(\theta^L+\theta^R)\cdot\xi
\end{equation}
Equation \eqref{eq:cond3} requires
$\ker\om_m\supseteq\{\xi_M(m)|\,(\on{Ad}_{\Phi(m)}+1)\xi=0\}$; 
for this reason the condition in Definition  \ref{def:qham}
is called the {\em minimal degeneracy 
condition}. If $G=T$ is a torus, one
recovers the usual definition of a symplectic $T$-manifold with torus
valued moment map.

\subsection{Conjugacy classes, symmetric spaces}
Basic examples for spaces with group valued moment maps
are $G$-conjugacy classes $\Co$. The moment map $\Phi$
is the embedding into $G$ and the 2-form is given by
\begin{equation}\label{eq:qKKS}
\om(\xi_\Co(g),\zeta_\Co(g))= -\hh\ (\Ad_{g}-\Ad_{g^{-1}})\xi\cdot\zeta 
\end{equation}
($g\in\Co$). 
The homogeneous group valued Hamiltonian $G$-manifolds are exactly 
the $G$-equivariant covering spaces of conjugacy classes.
\label{subsec:symmetric}
More generally, given a Lie group automorphism $\psi\in\on{Aut}(G)$
one defines {\em twisted conjugacy classes} to be the orbits 
of the action 
$$ \Ad^\psi_h(g)=\psi(h)gh^{-1}.$$
The twisted conjugacy classes are $\psi$-invariant 
since $\psi(g)=\Ad_g^\psi(g)$. 
Note $\Ad^\psi_G(e)=G/G^\psi$. Taking $\psi$ to be 
an involution we see that every compact symmetric space 
for $G$ is a twisted conjugacy class, up to finite cover.

Suppose $\psi$ has finite order $k$, and embed $\Z_k\hra \on{Aut}(G)$
as the subgroup generated by $\psi$.  For any $g\in G$, the map $G\to
\Z_k\ltimes G,\ g\mapsto (\psi^{-1},g)$ takes the twisted conjugacy
class of $G$ to conjugacy classes of the disconnected group
$\Z_k\ltimes G$.  In particular, if $\psi=\psi^{-1}$ is an 
involution, the identification of the symmetric space $G/G^\psi$
as a conjugacy class of $\Z_2\ltimes G$ is given by the map 
$\Phi:\,xG^\psi\mapsto (\psi,\psi(x)x^{-1})$. 
Note that the 2-form $\om$ given by \eqref{eq:qKKS} vanishes
for symmetric spaces.  The group $G$ itself is a symmetric space for
$G\times G$, where $\psi(a,b)=(b,a)$. Thus $G$ becomes a group valued
Hamiltonian $\Z_2\ltimes (G\times G)$ space, with 2-form $\om=0$,
moment map $a\mapsto (\psi,a,a^{-1})$ and action $(g_1,g_2).a=g_2 a
g_1^{-1},\ \ \psi.a=a^{-1}$.
\subsection{The four-sphere}
Let $G=\on{SU}(2)$ and $M=S^4$ the unit sphere in 
$\R^5\cong \C^2\times\R$, with $\on{SU}(2)$-action induced from 
the defining action on $\C^2$. In Appendix A we show that 
$M$ carries the structure of a group valued Hamiltonian 
$\on{SU}(2)$-manifold, with moment map $\Phi:\,S^4\to 
\on{SU}(2)\cong S^3$ the suspension of the Hopf fibration 
$S^3\to S^2$. 
\subsection{Exponentials}\label{subsec:exp}
Let $(M,\om_0,\Phi_0)$ be a Hamiltonian $G$-manifold in the usual sense. 
That is, $\om_0$ is a $G$-invariant symplectic form, and
$\Phi_0:\,M\to \g^*$ is an equivariant map satisfying the moment map 
condition, $\iota(\xi_M)\om_0=\d\Phi_0\cdot\xi$, using 
the inner product to identify $\g^*\cong\g$. Let
$\varpi\in\Om^2(\g)$ be the image of $\exp^*\eta$ under the de Rham 
homotopy operator $\Om^\star(\g)\to\Om^{\star-1}(\g)$. The triple 
$(M,\om_0+\Phi_0^*\varpi,\exp(\Phi_0))$ satisfies the 
axioms for a group valued 
Hamiltonian $G$-manifold, except possibly for the minimal degeneracy 
condition. It turns out that
$\om=\om_0+\Phi_0^*\varpi$ is minimally degenerate if and only if 
the exponential map has maximal rank at all points of $\Phi_0(M)$.
This construction takes (co-)adjoint $G$-orbits in $\g^*\cong\g$
to conjugacy classes for $G$. 

\subsection{Products}\label{sec:fusion}
Suppose 
$(M,\om,(\Phi_1,\Phi_2))$
is a group valued Hamiltonian $G\times G$-manifold. Then 
$\ti{M}=M$ with diagonal action, moment map $\ti{\Phi}=\Phi_1\Phi_2$, 
and 2-form 
$$ \ti{\om}=\om-\hh \Phi_1^*\theta^L\cdot\Phi_2^*\theta^R$$
is a group valued Hamiltonian $G$-manifold. If $M=M_1\times
M_2$ is a direct product of two group valued Hamiltonian $G$-manifolds, 
we call $\ti{M}=M_1\fus M_2$ the fusion product of $M_1,M_2$. 
\subsection{The double}
View $G$ as a group valued Hamiltonian $\Z_2\ltimes(G\times G)$-manifold
as above. The fusion product $G\fus G$ has moment map 
$$ (a_1,a_2)\mapsto (\psi,a_1,a_1^{-1})(\psi,a_2,a_2^{-1})
=(e,a_1^{-1}a_2,a_1a_2^{-1}).
$$ 
Since the moment map takes values in $G\times G\subset
\Z_2\ltimes(G\times G)$, we can restrict 
to the $G\times G$-action, and consider $G\fus G$ as a group valued Hamiltonian 
$G\times G$-manifold. Denoting $a:=a_1^{-1},b:=a_2$ the moment map reads
$$ (a,b)\mapsto (ab,a^{-1}b^{-1})$$
and the $G\times G$-action is
$$ (g_1,g_2).(a,b)=(g_1 a g_2^{-1}, g_2 b g_1^{-1}).$$
This is the {\em double} $D(G)$ introduced in \cite{al:mom}. 
Passing to the diagonal action, we obtain a group valued 
Hamiltonian $G$-manifold $\ti{D}(G)$, with $G$ acting 
by conjugation on each factor and moment map the Lie group 
commutator, $(a,b)\mapsto [a,b]=aba^{-1}b^{-1}$.

\subsection{Reduction} Suppose $(M,\om,(\Phi,\Psi))$ 
is a group valued Hamiltonian $G$-manifold. Then $e\in G$ is a regular
value of $\Phi$ if and only if the action of $G$ on $\Phi^{-1}(e)$ is
locally free, and in this case the {\em symplectic quotient} $M\qu
G=\Phi^{-1}(e)/G$ is a symplectic orbifold, with symplectic form
induced from $\om$. If one drops the regularity assumption, $M\qu G$
acquires more serious singularities and is a stratified symplectic
space in the sense of Sjamaar-Lerman. More generally, if $\Co\subset
G$ is a conjugacy class, let $\Co^-$ be its image under the inversion
map $g\mapsto g^{-1}$.  Then $e$ is a regular value for the action on
$M\fus \Co^-$ if and only if $\Co$ is contained in the set of regular
values of $\Phi$. We define $M_\Co=(M\fus \Co^-)\qu G \cong
\Phinv(\Co)/G$.

\section{Liouville measures}
Since the 2-form $\om$ of a group valued Hamiltonian $G$-manifold
$(M,\om,\Phi)$ may be degenerate, its top exterior power does not
define a volume form, in general. We will show in this Section that if
the adjoint action $G\to \on{O}(\g)$ lifts to the group $\on{Pin}(\g)$
(e.g. if $G$ is 1-connected), then $M$ carries a volume form $(e^\om
\Phi^*\ca{T})_{[top]}$ where $\ca{T}$ is a certain differential form
on the group $G$.  In the absence of such a lift, one still obtains a
nowhere vanishing measure $\nu$ on $M$, called the Liouville measure.

\subsection{Volume forms}
We recall the definition of the double covering 
$\on{Pin}(\g)\to \on{O}(\g)$ as a subset of the Clifford algebra
\cite{at:cli}. Let $\on{Cl}(\g)$ be the Clifford algebra of $\g$ with
respect to $-1/2$ the inner product on $\g$. If $e_i$ is a given 
orthonormal basis of $\g$, and $x_i$ the corresponding generators of
$\Cl(\g)$, the defining relations of $\Cl(\g)$ read, 
$ x_ix_j+x_jx_i=\delta_{ij}$.
The group $\Cl(\g)^\times$ of invertible elements acts on $\Cl(\g)$ 
by the twisted adjoint action  
$$ \wt{\Ad}(x)y=\alpha(x)yx^{-1},\ \ x\in \Cl(\g)^\times,\ y\in\Cl(\g)$$
where $\alpha:\,\Cl(\g)\to \Cl(\g)$ is the parity operator, i.e. 
$\alpha(x)=x$ for $x$ even and $\alpha(x)=-x$ for $x$ odd.
Let $\Pin(\g)\subset \Cl(\g)^\times$ be the subgroup 
generated by elements $\xi\in\g\subset \Cl(\g)$ with $\xi\cdot\xi=2$.
The action $\wt{\Ad}$ restricts to an 
isometry of $\on{Pin}(\g)$ on $\g\subset\Cl(\g)$, and 
$\wt{\Ad}:\,\on{Pin}(\g)\to \on{O}(\g)$ is 
a double covering. The group $\Pin(\g)$ has two connected
components, given as intersections
$\Spin(\g)=\Pin(\g)\cap 
\Cl(\g)_{\on{even}}$ and 
$\Pin(\g)\cap \Cl(\g)_{\on{odd}}$.

Let us assume that the adjoint action $\Ad:\ G \to \on{O}(\g)$ 
lifts to a group homomorphism $\tau: G \to \Pin(\g)$:
%
$$\xymatrix{&\Pin(\g)\ar[d]\ar[r]&\Cl(\g)\\ 
G \ar[r]^{\Ad} \ar@{.>}[ru]^{\tau}& \on{O}(\g)&}$$
The restriction of $\tau$ to the identity component takes values 
in $\on{Spin}(\g)$ 
and is given by the formula 
\begin{equation}\label{DefTau}
\tau(\exp\,\mu)=\exp(-\hh \sum_{i,j,k} f_{ijk}\,\mu_i\,x_j\,x_k).
\end{equation}
where $f_{ijk}=[e_i,e_j]\cdot e_k$ are the structure constants.
Let $\wedge\g$ be the exterior algebra, with generators $y_i$
and relations $y_iy_j+y_jy_i=0$, 
and $\sig:\,\on{Cl}(\g)\to\wedge\g$ the symbol map, 
taking $x_{i_1}\cdots x_{i_s}$ to 
$y_{i_1}\cdots y_{i_s}$ for $i_1<\cdots<i_s$. 
Denote by $s_{1/2}:\,\wedge\g\to\wedge\g$ the linear map 
equal to multiplication by $2^{-l}$ on $\wedge^l\g$. We define 
$\ca{T}\in\Om(G)$ to be the differential form 
\begin{equation}\label{eq:T}
\ca{T}:=s_{1/2}\circ \sig\circ \tau \in C^\infty(G,\wedge\g)^G
\cong \Om(G)^G,
\end{equation}
using the isomorphism defined by right-invariant (or equivalently 
left-invariant) differential forms. That is, the value of $\ca{T}$ 
at $g$ is obtained from $\sig(\tau(g))$ by replacing the variables 
$y_i$ with $\hh (\theta_i^R)_g$, where $\theta^R_i=\theta^R\cdot e_i$
are the components of the right-invariant Maurer-Cartan forms. 
Our main theorem reads:
\begin{theorem}\label{th:VolumeForm}
For any group valued Hamiltonian $G$-manifold $(M,\om,\Phi)$, 
the top form degree part of the differential form 
\begin{equation}\label{eq:defi}
 \Gamma=\exp(\om)\Phi^*\ca{T}.
\end{equation}
is a volume form on $M$. On the subset where $\om$ is non-degenerate, 
we have
\begin{equation}\label{eq:formula1}
(\Gamma)_{[\dim M]}=\pm 
\f{\exp(\om)_{[\dim M]}}
{\big|{\det}(\f{1+\Ad_{\Phi}}{2})\big|^{1/2} }.
\end{equation}
The signs depend on the choice of lift $\tau$. 
\end{theorem}
The proof of this Theorem is given in Section \ref{sec:proofs}.
\begin{remark}
The volume form 
depends not only on the 2-form but also on the moment map. 
If $c\in Z(G)$ is a central element, 
then $\Phi'=c\Phi$ is a moment map for the 
same group action and 2-form $\om$.
Since $\tau$ is a group homomorphism, $\Gamma'=\tau(c)\Gamma$. 
We have $\tau(c)=\pm 1$ since $\Ad_c=1$.
Therefore, the new volume form $\Gamma'_{[\dim M]}$ differs from 
$\Gamma_{[\dim M]}$  by the sign of $\tau(c)$.
\end{remark}
If $G$ is connected, we can be more precise. Let $T$ be a maximal 
torus of $G$, with Lie algebra $\t$. 
Choose a system of positive (real) roots, and let 
$\rho\in\t^*$ be their half-sum. 
It is known that the adjoint action $\Ad:\ G\to \on{SO}(\g)$ 
lifts to a homomorphism $\tau:\ G\to\Pin(\g)$ if and only if 
$\rho$ is in the lattice of weights for $T$.
\footnote{In fact, $\rho$ is a weight if and only if 
$G/T$ admits a $G$-equivariant Spin structure (cf. 
\cite{goe:eq} for discussion and references).} 
The character $\chi_\rho$ of the $\rho$-representation
$V_\rho$ defines a smooth square root of $g\mapsto
\det(\Ad_g+1)$:
$$ \f{\chi_\rho(g)}{\dim V_\rho}= {\det}^{1/2}(\f{\Ad_g+1}{2}).$$
We have the following extension of Theorem \ref{th:VolumeForm}:
\begin{theorem}\label{cor:conn}
Suppose $G$ is a compact connected group, and assume that $\rho$ 
is a weight for $G$. Let $(M,\om,\Phi)$ be a 
group valued Hamiltonian 
$G$-manifold $(M,\om,\Phi)$. Then $M$ is even-dimensional, and 
carries a canonical volume form. 
On the subset where  $\Phi^*\chi_\rho\not=0$, the volume form is 
given by
$$ 
\f{\dim V_\rho}{\Phi^*\chi_\rho}\ \exp(\om)_{[\dim M]}.
$$
\end{theorem}

\subsection{Proof of Theorems \ref{th:VolumeForm}
and \ref{cor:conn}}\label{sec:proofs}
\label{subsec:cons}
We first give a more explicit description of the differential form $\ca{T}$.
For any subspace $\mf{s}\subset \g$ with a given orientation let
$$\d\vol_{\mf{s}}\in 
\wedge^{\dim{\mf{s}}}\mf{s}$$
denote the Riemannian volume form. Also, for any endomorphism $A$ of 
$\g$ preserving $\mf{s}$ we denote $\det_{\mf{s}}(A):=\det(A|_\mf{s})$.
Given $g\in G$ let  
$\g^{\pm}$ denote the $\pm 1$ eigenspaces of
$\Ad_g:\,\g\to\g$, 
and $\g'$ the orthogonal complement of $\g^{-}\oplus\g^{+}$. Thus 
\begin{equation}\label{eq:split}
\g=\g'\oplus \g^{-}\oplus \g^{+}.
\end{equation}
Accordingly, decompose $I=(1,\ldots,\dim\g)$ as $I=I'\cup I^{-}\cup
I^{+}$ where $I'$ denotes the first $\dim\g'$ indices, $I^{-}$ the
following $\dim\g^{-}$ indices, and $I^{+}$ the remaining $\dim\g^{+}$
indices. Choose an orthonormal basis $e_i,\,i\in I$ of $\g$ such that
the basis vectors labeled by $I',I^-,I^+$ span $\g',\g^-,\g^+$,
respectively.  For any linear map $A:\,\g\to\g$ we denote by
$A_{ij}=e_i\cdot A(e_j)$ its components. We denote by
$\theta^R_i=\theta^R\cdot e_i$ and $\theta^L_i=\theta^L\cdot e_i$ the
components of the Maurer-Cartan forms.  On $\g'$ the operator
$\f{\Ad_g-1}{\Ad_g+1}$ is well-defined and invertible.
\begin{lemma}\label{lem:dec}
At any point $g\in G$, the form $\ca{T}$ is a product 
$\pm \ca{T}'\ca{T}^{-}$, where 
\begin{equation}\label{eq:T1}
 \ca{T}^{-}= 2^{-\dim\g^{-}/2}\ \prod_{i\in I^-}\theta_i^R
\end{equation}
and 
\begin{equation}\label{eq:T2}
 \ca{T}'=\big| {\det}_{\g'}\big(\f{\Ad_g+1}{2}\big)\big|^{1/2}
\exp\Big(-\f{1}{4} \sum_{i,j
\in I'}\big(\f{\Ad_g-1}{\Ad_g+1}\big)_{ij}
\theta_i^R\theta_j^R\Big). 
\end{equation}
(The sign $\pm$ depends on the choice of lift and on the orientation
on $\g^{-}$ defined by the choice of basis.)
If $G$ is connected we have the formula, valid on 
the subset where $\chi_\rho\not=0$, 
\begin{equation}\label{eq:nice}
\ca{T}=\f{\chi_\rho}{\dim V_\rho}\ \exp\Big(-\f{1}{4} 
\sum_{i,j}\big(\f{\Ad_g-1}{\Ad_g+1}\big)_{ij}
\theta_i^R\theta_j^R\Big).
\end{equation}
\end{lemma}
\begin{proof}
Write $\Ad_g=:S=S' S^{-}$, where $S^{-}=-\on{Id}_{\g^-}\in
\on{O}(\g^{-})$ and $S'\in \on{SO}(\g')$. The product
$$\ti{S}^-=2^{\dim\g^-/2} \prod_{i\in I^-}x_i$$
is a lift of $S^-$, with symbol $\sig(\ti{S}^{-})=2^{\dim\g^{-}/2}
\prod_{i\in I^-} y_i$.  
As a special case of \cite[Proposition 3.13]{be:he}, 
the symbol of any lift $\ti{S}'\in \Spin(\g')$ of $S'$ is given by a formula 
$$ \sig(\ti{S}')=\pm\  \big|{\det}\big(\f{S'+1}{2}\big)\big|^{1/2}
\exp\Big(-\sum_{i,j\in I'}
\big(\f{S'-1}{S'+1}\big)_{ij}y_i y_j\Big),$$
where the sign depends on the choice of lift. 
Replacing the variables 
$y_i\in\g$ with $\hh (\theta_i^R)$, the first part of the 
Lemma follows. Equation \eqref{eq:nice} is a special case 
since $\g^{-}=\{0\}$ on the set where $\chi_\rho\not=0$; 
the sign is verified by evaluating at $g=e$. 
\end{proof}

\begin{proof}[Proof of Theorem \ref{th:VolumeForm}]
Given $m\in M$ let $g=\Phi(m)$.  From \eqref{eq:split}
we obtain the following splitting of the tangent space 
\begin{equation} \label{eq:firstsplit}
T_mM=E\oplus \g_g^\perp=E\oplus \g'\oplus\g^-. 
\end{equation}
Here $\g_g^\perp=\g'\oplus\g^-$ 
is the orthogonal complement of the stabilizer algebra $\g_g$, embedded 
by the generating vector fields $\xi\mapsto \xi_M(m)$,  
and 
$E=\{v\in T_mM|\,\iota(v)\Phi^*\theta^R\in \g_g\}$ 
is the pre-image under $\d_m\Phi$ of the tangent space of the stabilizer group 
$G_g$. The moment map 
condition shows that this splitting is $\om$-orthogonal and 
identifies
$\ker\om_m\cong\g_-$. In particular, the 2-forms 
$\om_E=\om|_E$ and $\om'=\om|_{\g'}$ are symplectic. 

We will show that in terms of the splitting 
$T_mM=E\oplus \g_g^\perp$, 
the value of $\Gamma$ at $m\in M$ is given by the formula 
\begin{equation}\label{eq:first}
(\Gamma_m)_{[\dim M]}=\pm 
|{\det}_{\g_{g}^\perp}(\Ad_{g}-1)|^{1/2}\ 
(\exp\om_E)_{[\dim E]}\wedge \d\vol_{\g_{g}^\perp}.
\end{equation}
In particular, \eqref{eq:first} shows 
that $\Gamma_{[\dim M]}$ is a volume form. 
Consider the splitting $T_mM=E\oplus\g'\oplus\g^-$ and 
the corresponding decomposition
$\Gamma=\pm e^{\om_E}\,\Gamma'  \Gamma^-$
with $\Gamma'=(\exp\om')\ca{T}'$
and $\Gamma^{-}=\ca{T}^{-}$. 
(We are dropping the base point $m\in M$ from the notation.)
The form $\om'$ 
is determined by the moment map condition, and is given by 
$$\om'=\f{1}{4} 
\sum_{i,j\in I'}
\big(\f{\Ad_g+1}{\Ad_g-1}\big)_{ij}\theta_i^R 
\theta_j^R.
$$
Indeed, since $\iota_i\theta^R_j=(\Ad_g-1)_{ij}$, it 
is easily verified that 
$ \iota_i\om'=\hh (\theta^L+\theta^R)_i$ for all 
$i\in I'$, as required.  
Using \eqref{eq:T2} and the following 
equality of operators on $\g'$, 
$$ 
\f{\Ad_g+1}{\Ad_g-1}-\f{\Ad_g-1}{\Ad_g+1}=\f{4}{\Ad_g-\Ad_{g^{-1}}},
$$ 
we obtain 
$$ \Gamma'=\pm 
\big|{\det}_{\g'}\ \big(\f{\Ad_g+1}{2}\big)\big|^{1/2}
\exp\big(\sum_{i,j
\in I'}
\big(\f{1}{\Ad_g-\Ad_{g^{-1}}}\big)_{ij}
\theta_i^R\theta_j^R \big).$$
Thus
\begin{eqnarray}\nonumber
(\Gamma')_{[\dim \g']}
&=&\pm \ \big|{\det}_{\g'}\ \big(\f{\Ad_g+1}{2}\big)\big|^{1/2}\,
\big|{\det}_{\g'}\big(\f{2}{\Ad_g-\Ad_{g^{-1}} }\big)\big|^{1/2}
\prod_{i \in I'}\theta_i^R\\
\label{eq:emma}
&=&\pm \ \big|{\det}_{\g'}\big({\Ad_g-1}\big)\big|^{-1/2}
\prod_{i \in I'}\theta_i^R.
\end{eqnarray}
This expression combines with $\Gamma^{-}=\ca{T}^{-} $ to a factor,
$$ (\Gamma'\Gamma^{-})_{[\dim\g_g^\perp]}
=\big|{\det}_{\g_g^\perp}\big({\Ad_g-1}\big)\big|^{-1/2}
\prod_{i \in I'\cup I^{-}}\theta_i^R
=\big|{\det}_{\g_g^\perp}\big({\Ad_g-1}\big)\big|^{1/2}
\d\vol_{\g_g^\perp},
$$
where we have used \eqref{eq:T1} with 
$2^{\dim \g^{-}}=|{\det}_{\g^{-}}(\Ad_g-1)|$, 
and $\iota_i\theta_j^R=(\Ad_g-1)_{ij}$.
This proves \eqref{eq:first}. For \eqref{eq:formula1}, assume that 
$\om_m$ is non-degenerate, that is $\ker\om_m\cong \g^{-}=\{0\}$. 
We have
$$ (\exp\om')_{[\dim \g']}
=\pm \ \big|{\det}_{\g'}\big(\f{1}{2}\
\f{\Ad_g+1}{\Ad_g-1}\big)\big|^{1/2}
\prod_{i \in I'}\theta_i^R.
$$
Comparing with \eqref{eq:emma}, this shows 
$$(\Gamma')_{[\dim \g']} 
=\pm \ 
{(\exp \om')_{[\dim \g']}}\
{\big|{\det}_{\g'}\big(\f{\Ad_g+1}{2}\big)\big|^{-1/2}}
$$
which yields \eqref{eq:formula1}. 
\end{proof}

\begin{proof}[Proof of Theorem \ref{cor:conn}]
Since $G$ is connected, the lift $\tau:\,G\to \on{Spin}(\g)$ is unique. 
Moreover, since $\on{Spin}(\g)\subset \Cl(\g)_{\on{even}}$, the 
form $\ca{T}$, hence also $\Gamma$, vanishes in odd degrees. 
The formula for $\Gamma_{[\dim M]}$ is obtained by 
using the expression \eqref{eq:nice} for $\ca{T}$ in
the proof of Theorem \ref{th:VolumeForm}, and keeping track of the signs.
\end{proof}

\subsection{Products}
\, From the differential form expression \eqref{eq:defi}, it is not
obvious how the volume forms behave under the fusion operation from
Section \ref{sec:fusion}. In particular, there is no simple relation
of the product $\Phi_1^*\ca{T}\,\Phi_2^*\ca{T}$ with
$(\Phi_1\Phi_2)^*\ca{T}$.  It will be convenient to work with an
equivalent expression for $\ca{T}$ involving the Clifford algebra
$\Cl(\g)$.
Let $P_{\on{hor}}^{\on{Cl}}:\,\Cl(\g)\to\R$ denote
horizontal projection, i.e. the linear map defined by 
$P_{\on{hor}}^{\on{Cl}}(1)=1$ and 
$P_{\on{hor}}^{\on{Cl}}(x_{j_1}\ldots x_{j_s})=0$,
where $j_1 < j_2 < \cdots < j_s$.
For any $\Z_2$-graded space $\A$ it extends to a linear map
$P_{\on{hor}}^{\on{Cl}}:\ \Cl(\g)\otimes \A\to \A$ (graded 
tensor product). 
Observe  that for any $j_1 < j_2 < \ldots < j_s$,  the operator
$P_{\on{hor}}^{\on{Cl}}:\, \Cl(\g)\otimes\wedge\g\to\wedge\g$ 
takes 
$$e^{-2\sum_i x_i y_i}x_{j_1}\cdots x_{j_s}
=\prod_i (1-2 x_i y_i) x_{j_1}\cdots x_{j_s}
=(x_{j_1}+y_{j_1})\cdots (x_{j_s}+y_{j_s})
$$
to $y_{j_1}\cdots y_{j_s}=\sig(x_{j_1}\cdots
x_{j_s})$. This shows 
$P_{\on{hor}}^{\on{Cl}}(e^{-2\sum_i x_i y_i}\tau)=\sig(\tau)$.
Replacing $y_i$ with $\hh\theta_i^R$ we obtain, 
\begin{equation}\label{eq:nicer}
 \ca{T}=P_{\on{hor}}^{\on{Cl}}(e^{-\sum_i x_i\theta_i^R}\tau).
\end{equation}
\begin{theorem}[Products]\label{th:fusion}
Let $(M,\om,(\Phi_1,\Phi_2))$ be a group valued Hamiltonian $G$-manifold 
and $(\ti{M},\ti{\om},\ti{\Phi})$ its fusion as defined in \ref{sec:fusion}.
The corresponding forms $\Gamma,\ti{\Gamma}$ are related by 
\begin{equation}\label{GG}
\ti{\Gamma}=\exp(-\hh \sum_i \iota_i^{1}\iota_i^{2})\,\Gamma.
\end{equation}
where $\iota_i^1$ are the contraction operators for the action of the 
first $G$-factor on $M$, and $\iota_i^2$ those for the second $G$-factor. 
In particular $\ti{\Gamma}_{[\dim M]}=\Gamma_{[\dim M]}$. 
\end{theorem}
It follows that the volume form of a fusion product $M_1\fus M_2$ is 
the product of the volume forms of $M_1,M_2$.
\begin{proof}
We will need the following two identities from 
\cite[Section 5.3]{al:no}, both of which are verified by 
straightforward calculation:
\begin{equation}\label{eq:11}
\Phi_1^*(e^{-\sum_i x_i\theta_i^{R}}\tau) \Phi_2^*(e^{-\sum_i x_i\theta_i^{R}}\tau)
=e^{-\hh \Phi_1^*\theta^L \cdot \Phi_2^*\theta^R}
(\Phi_1\Phi_2)^*(e^{-\sum_i x_i\theta_i^{R}}\tau)
\end{equation}
\begin{equation}\label{eq:12}
 (\iota_i^{\Cl}+\iota_i+ \hh(\theta_i^L+\theta_i^R))
\big(e^{-\sum_j x_j\theta_j^{R}}\tau\big)=0.
\end{equation}
Here $\iota_i:\,\Om(G)\to\Om(G)$ are the contraction operators for
the conjugation action, and $\iota_i^{\Cl}=\ad(x_i)$ the
contraction operators for the Clifford algebra.  
From \eqref{eq:12} and the moment map condition we obtain 
\begin{equation}\label{eq:contr}
(\iota_i^{\Cl}+\iota_i)
\big(e^{\om}\Phi^*(e^{-\sum_j x_j\theta_j^{R}}\tau)\big)=0.
\end{equation}
Since $\ti{\om}=\om-\hh \Phi_1^*\theta^L \cdot \Phi_2^*\theta^R$, 
Equation \eqref{eq:11} gives 
\begin{eqnarray}\label{eq:21}
\ti{\Gamma}&=&
P_{\on{hor}}^{\Cl}
\big(e^\om\Phi_1^*(e^{-\sum_i x_i\theta_i^{R}}\tau) 
           \Phi_2^*(e^{-\sum_i x_i\theta_i^{R}}\tau)\big)\\
\nonumber 
&=&
P_{\on{hor}}^{\wedge\g}
\sig\big(e^\om\Phi_1^*(e^{-\sum_i x_i\theta_i^{R}}\tau) 
           \Phi_2^*(e^{-\sum_i x_i\theta_i^{R}}\tau)\big)
\end{eqnarray}
where $P_{\on{hor}}^{\wedge\g}$ is horizontal projection for the 
exterior algebra $\wedge\g$. Let $\iota_i^{\wedge\g}=\iota(e_i)$ be the 
contraction operators for the exterior algebra.
The symbol map takes 
the Clifford product on $\Cl(\g)$ to the product on 
on $\wedge\g$ given 
by application of the operator
$\exp(-\hh \sum_i \iota_i^{\wedge\g,1}\iota_i^{\wedge\g,2})$
on $\wedge\g\otimes\wedge\g$ (where the superscripts refer
to the $\wedge\g$-factor), followed by
wedge product (see \cite[Theorem 16]{ko:cl}, or \cite[Lemma
3.1]{al:no}). Because of \eqref{eq:contr} we may 
replace $\exp(-\hh \sum_i \iota_i^{\wedge\g,1}\iota_i^{\wedge\g,2})$
with $\exp(-\hh \sum_i \iota_i^{1}\iota_i^{2})$ everywhere, 
which then commutes with $P_{\on{hor}}^{\wedge\g}$. Equation 
\eqref{eq:21} becomes
$$ \ti{\Gamma}=e^{-\hh \sum_i \iota_i^{1}\iota_i^{2}}
P^{\wedge\g}_{\on{hor}}
\big(e^\om \sig(\Phi_1^*(e^{-\sum_i x_i\theta_i^{R}}\tau)) 
           \sig(\Phi_2^*(e^{-\sum_i x_i\theta_i^{R}}\tau))\big).
$$
Since $P_{\on{hor}}^{\wedge\g}$ is a ring homomorphism, 
this is equal to 
$\exp(-\hh \sum_i \iota_i^{1}\iota_i^{2})(e^\om\Phi_1^*\ca{T}
\Phi_2^*\ca{T})$, proving the Theorem.
\end{proof}

\subsection{General case}\label{sec:general}
For general compact Lie groups $G$, a lift $\tau:\,G\to \on{Pin}(\g)$
of the adjoint action need not exist. We will now show that for any
group valued Hamiltonian $G$-manifold $(M,\om,\Phi)$, the form
$\Gamma$ may be defined as a form with values in the orientation
bundle $o_M$.  Recall that $o_M$ is the real line bundle associated to
the oriented double cover of $M$. One defines the space of twisted
differential forms $ \Om_t(M)=\Om(M,o_M)$ as sections of $\wedge
T^*M\otimes o_M$.  The space $\Om_t(M)$ is naturally isomorphic to the
space of differential forms on the oriented double cover of $M$ that
are anti-invariant under deck transformations. The real line bundle
$\wedge^{\dim M} T^*M \otimes o_M$ is isomorphic to the density bundle
of $M$; hence the space $\Om^{\dim M}_t(M)$ is just the space of
densities (smooth measures) of $M$.

Define double coverings $\pi_G:\,\wh{G}\to G$ and $\pi_M:\,\wh{M}\to M$ 
by the pull-back diagram
$$
\xymatrix{\wh{M}\ar[d]^{\pi_M}
\ar[r]^{\hat{\Phi}}&\wh{G}\ar[d]^{\pi_G}\ar[r]^{\hat{\tau}}&
\Pin(\g)\ar[d]^\pi\\ 
M\ar[r]^{{\Phi}}& G \ar[r]^{\Ad} \ar[r]^{\Ad}& \on{O}(\g)&}$$
Write $\ker\pi_G=\{e,c\}$.  The conjugation action of $G$ on $\wh{G}$
and the action of $G$ on $M$ define an action on the fiber product
$\wh{M}$.  Let $\wh{\om}$ be the pull-back of $\om$. Then
$(\wh{M},\wh{\om},\wh{\Phi})$ is a group valued Hamiltonian $\wh{G}$
space, with $\ker\pi_G$ acting trivially. The construction in
\ref{subsec:cons} defines a differential form $\wh{\Gamma}$ on
$\wh{M}$, for which $\wh{\Gamma}_{[\dim M]}$ is a volume form.

Let $S:\,\wh{M}\to\wh{M}$ be the non-trivial 
deck transformation. Then $S$ is $\wh{G}$-equivariant, preserves the 2-form 
$\wh{\om}$, and changes the moment map by 
$S^*\wh{\Phi}=c \wh{\Phi}$. Thus 
$$S^*\wh{\Gamma}=\wh{\tau}(c)\wh{\Gamma}=-\wh{\Gamma}.$$
In particular, 
$S$ changes the orientation of $\wh{M}$. This identifies $\wh{M}$ 
with the oriented double cover of $M$. The form   
$\wh{\Gamma}$ descends to a 
form $\Gamma\in\Om_t(M)$, with top degree part 
a strictly positive measure. We call
$\nu:=\Gamma_{[\dim M]}$ the {\em Liouville measure}, and its 
integral $\Vol(M):=\int_M\nu$ the {\em Liouville volume}. 

The fusion formula \ref{th:fusion}
holds in this more general context; in particular, the Liouville 
measure $\ti{\nu}$ for the ``fused'' space $\ti{M}$ coincides with 
the original Liouville measure $\nu$. 

\subsection{Conjugacy classes} 
Going back to the proof of Theorem \ref{th:VolumeForm}, we 
obtain the following formula for the Liouville measure of a conjugacy 
class in $G$:
\begin{proposition}\label{prop:LioConj}
The Liouville volume of the conjugacy class
$\Co\subset G$ of $g\in G$ is given by the formula
\begin{equation}\label{eq:Eq42a}
\Vol(\Co)=|{\det}_{\g_g^\perp}(\Ad_g-1)|^\hh \f{\vol G}{\vol G_g}. 
\end{equation}
where $\vol G$ and $\vol G_g$ are the Riemannian volumes for the given
inner product on $\g$, and $\g_g^\perp$ is the orthogonal complement
of the stabilizer algebra of $g$.
\end{proposition} 
\begin{proof}
By passing to a cover if necessary, we may assume that $G$ admits a
lift $\tau:\,G\to \on{Spin}(\g)$. The formula follows directly from
Equation \eqref{eq:first} in the proof of Theorem \ref{th:VolumeForm},
since the subspace $E$ is trivial in this case.
\end{proof}

The formula \eqref{eq:Eq42a} is similar to the well-known formula for
the symplectic volume of a coadjoint orbit
$\O=G\cdot\mu\subset\g^*\cong\g$:
\begin{equation}\label{eq:Eq42b}
\Vol(\O)=|{\det}_{\g_\mu^\perp}
(\ad_\mu)|^\hh \f{\vol G }{\vol G_\mu}
\end{equation}

\subsection{Exponentials}
Let $J:\,\g\to\R$ the determinant of the Jacobian of the exponential 
map $\exp:\,\g\to G$, 
$$ 
J(\xi)={\det}_{\g_\xi^\perp}\big(\f{e^{\ad_\xi}-1}{\ad_\xi}\big),$$
and $J^{1/2}$ the unique smooth square root equal to $1$ at the
origin.  Suppose $J(\xi)\not=0$ so that $\exp$ has maximal rank at
$\xi$.  The conjugacy class $G.\exp(\xi)$ is the group valued
Hamiltonian $G$-manifold corresponding (in the sense of
\ref{subsec:exp}) to the (co)-adjoint orbit $G.\xi$. Equations
\eqref{eq:Eq42a} and \eqref{eq:Eq42b} show that
$$\f{\Vol(G.\xi)}{\Vol(G.\exp\xi)}=|J(\xi)^\hh|.$$
This generalizes, as follows:
\begin{proposition}\label{prop:exp}
Let $(M,\om_0,\Phi_0)$ be a Hamiltonian $G$-manifold, with 
Liouville measure $\nu_0$. Suppose
$\Phi_0^*J\not=0$ everywhere, and let $(M,\om,\Phi)$ the corresponding 
group valued Hamiltonian $G$-manifold. Then the Liouville measure of 
$(M,\om,\Phi)$ is given by 
$$ \nu=\Phi_0^*|J^{1/2}| \nu_0.$$
\end{proposition}
\begin{proof}
We may assume that $G$ is connected and that it admits 
a lift $\tau:\,G\to \on{Spin}(\g)$. We compare the two measures
at some point $m\in M$. 
Again we go back to 
Equation \eqref{eq:first} from the proof of Theorem 
\ref{th:VolumeForm}. Let $\xi=\Phi_0(m)$ 
and $\h=\g_{\xi}=\g_{\exp\xi}$.
The splitting $T_mM=E\oplus \h^\perp$ is both $\om_0$-orthogonal and 
$\om$-orthogonal. By the moment map condition, 
the restriction of $\om_0$ to $\h^\perp$ 
is given by the skew-adjoint operator $\ad_\xi$. 
Hence, its
top exterior power is  
$\pm|{\det}_{\h^\perp}(\ad_\xi)|^\hh\, \d\vol_{\h^\perp}$. 
Together with \eqref{eq:first} this proves the Proposition.
\end{proof}


\section{Duistermaat-Heckman measures}
\subsection{Volumes of symplectic quotients}\label{subsec:volumes}
For any compact group valued Hamiltonian $G$-manifold $(M,\om,\Phi)$
we define the {\em Duistermaat-Heckman (DH) measure} as the
push-forward of the Liouville measure $\nu=\Gamma_{[\dim M]}$ under
the moment map,
$$ \m=\Phi_*\nu \in\ca{E}'(G)^G.$$ 
By general properties of push-forwards, the singular 
support of $\m$ is contained in the set of singular 
values of $\Phi$. 
Similar to the DH-measure for $\g^*$-valued moment maps, the 
measure $\m$ encodes volumes of reduced spaces. 
Recall that $G.m\subset M$ is called a {\em principal orbit} 
if the stabilizer groups for all orbits in a neighborhood 
of $G.m$ are conjugate to $G_m$. 
The stabilizer group $H=G_m$ is called a {\em principal stabilizer}. 
The union $M_{\on{prin}}$ 
of principal orbits is open and dense in $M$.
If $M/G$ is connected, then 
$M_{\on{prin}}/G$ is connected and any two principal 
stabilizers are conjugate. 
\begin{theorem}\label{th:meaningofdh}
Let $(M,\om,\Phi)$ be a compact group valued Hamiltonian $G$-manifold.  
Assume that the orbit space 
$M/G$ is connected and that the principal stabilizer 
groups are finite. 
Suppose each component of $\Phinv(e)$ meets $M_{\on{prin}}$. 
Then $\m$ is continuous on a neighborhood of $e$, and 
the volume of the (possibly singular) symplectic quotient 
$M\qu G$ is given by
\begin{equation}\label{eq:formula}
 \Vol(M\qu G)=\f{k}{\vol G}\f{\m}{\d\vol_G}\Big|_e ,
\end{equation}
where $k$ is the cardinality of a principal stabilizer. 
\end{theorem}
\begin{proof}
We may assume that $(M,\om,\Phi)$ is the exponential of a 
Hamiltonian space $(M_0,\om_0,\Phi_0)$.  The symplectic quotients 
$M\qu G$ and $M_0\qu G$ coincide. 
Since $\nu=\Phi_0^*|J^{1/2}|\,\nu_0$, the DH-measures are related by 
$\m=|J^{1/2}|\,\exp_*\m_0$. Since $J(0)=1$, we have  
$\f{\m}{d\vol_G}\big|_e=\f{\m_0}{d\vol_\g}\big|_0$. Hence, 
\eqref{eq:formula} follows from the well-known statement for 
Hamiltonian $G$-manifolds. 
\end{proof}

The right hand side of \eqref{eq:formula} can be re-expressed as a 
Fourier series. Let $\on{Irr}(G)$ denote the set of equivalence classes 
of irreducible unitary $G$-representations. For $\lambda\in \on{Irr}(G)$ 
let $\chi_\lambda$ denote the character of the corresponding 
irreducible representation $V_\lam$. The measure $\m$ 
is determined by its Fourier coefficients $\l\m,\chi_\lam\r$, 
\begin{equation}\label{eq:FourierDecomp1}
\f{\dh}{\d\vol_G} =\f{1}{\vol G}\sum_{\lambda\in\on{Irr}(G)} \l\dh,
\chi_\lambda\r \chi_\lambda^*\  
\end{equation}
In general, this sum need not converge as a continuous
function at $g=e$. We therefore 
proceed as in Liu's paper \cite{li:he} and
apply  a smoothing operator $e^{t\Delta},\ t>0$ to $\dh$, 
where $\Delta$ is the Laplace-Beltrami operator on $G$. 
Then $\lim_{t\to 0^+}(e^{t\Delta}\m)|_e=\m|_e$, so that 
\begin{equation}\label{eq:FourierDecomp}
\Vol(M\qu G)
=\f{k}{\vol G^2}\lim_{t\to 0^+}\sum_{\lambda\in\on{Irr}(G)} 
e^{-t p(\lam)}\l\dh,\chi_\lambda\r \dim V_\lam  
\end{equation}
where $p(\lam)$ is the eigenvalue of $-\Del$ on the character  
$\chi_\lam$. If $G$ is connected, the irreducible representations are labeled 
by the set of dominant (real) weights $\lambda\in \t^*$, and 
$$p(\lam)=\Vert\lambda+\rho\Vert^2-\Vert\rho\Vert^2.$$
In the general case, $V_\lam$ 
splits into a direct sum of irreducible representations of
the identity component $G^{\on{o}}$ with highest
weights $\lambda_i$, and $p(\lambda)= \Vert\lambda_i+\rho\Vert^2-
\Vert\rho\Vert^2$ for any $i$. 

If the Fourier coefficients $\l\m,\chi_\lam\r$ decrease sufficiently 
fast as $p(\lambda)\to \infty$, the convergence factor may be omitted 
and one has a simpler formula 
\begin{equation}\label{eq:volform1}
\Vol(M\qu G)=\f{k}{\vol G^2}\sum_{\lambda\in\on{Irr}(G)} \l\m,\chi_\lambda\r 
\dim V_\lambda.
\end{equation}
Suppose $G$ is connected, and view $\lambda$ as a weight. By the 
Weyl dimension formula, $\dim V_\lam$ is a polynomial of degree 
${\dim G-\on{rank} G}$ in $\lambda$.
Hence a sufficient criterion for absolute convergence is 
$\l\m,\chi_\lambda\r \le C ||\lambda||^{-\dim G-\eps}$ for some constants 
$C>0,\eps>0$.

\subsection{Fusion}
By Theorem \ref{th:fusion}, the Liouville measures $\nu$ and $\ti{\nu}$
on a group valued Hamiltonian $G\times G$ space $M$ and
on its fusion $\ti{M}$  coincide. Therefore, if $M$ is compact, 
the DH-measures $\m\in\ca{E}'(G\times G)$ and $\ti{\m}\in\ca{E}'(G)$ 
are related by push-forward under 
group multiplication $\Mult_G:\,G\times G\to G$, and the 
Fourier coefficients are
\begin{equation}\label{eq:conv0}
\l\ti{\m},\chi_\lam\r=(\dim V_\lam)^{-1} \l\m,\chi_\lam\otimes\chi_\lam\r.
\end{equation}
In particular, 
the DH-measure for a fusion product $M_1\fus M_2\cdots \fus M_r$ is 
given by convolution, $\m_1 * \m_2 *\cdots *\m_r$, with
Fourier coefficients 
\begin{equation}\label{eq:conv}
 \l\m_1 * \cdots *\m_l,\chi_\lam\r=
(\dim V_\lam)^{1-l}\prod_{j=1}^l \l\m_j,\chi_\lam\r.
\end{equation}
\subsection{Examples}

In this section we compute the DH-measures for a number of examples,
including the space $G^{2h}\times\Co_1\times\cdots\times\Co_r$ 
related to the moduli space of
flat bundles on a surface of genus $h$.  We begin with twisted
conjugacy classes.
\begin{proposition}
Let $\psi\in\on{Aut}(G)$ be an 
automorphism of finite order $k$, and use the same notation for 
the induced automorphism of $\g$. Given $g\in G$, let 
$H$ be the stabilizer of $g$ under the twisted adjoint action. 
The Liouville measure of the twisted 
conjugacy class $\Co=\Ad^\psi_G(g)$ is equal to the Riemannian 
measure of the homogeneous space $G/H$ times a factor 
$|{\det}_\mf{c}(\Ad_g-\psi)|^\hh$, where $\mf{c}\subset\g$ 
is the orthogonal complement of the kernel of 
$(\Ad_g-\psi)$. That is, 
$$ \Vol(\Co)= |{\det}_\mf{c}(\Ad_g-\psi)|^\hh \f{\vol G}{\vol H}. 
$$
\end{proposition}
\begin{proof}
View $\Co$ as a conjugacy class of $(\psi^{-1},g)$ 
for the semi-direct product $\Z_k\ltimes G$ defined by $\psi$.
The adjoint action of $(\psi^{-1},g)$ on $\g$ is given by
$\Ad_{(\psi^{-1},g)}\xi=\psi^{-1}(\Ad_g(\xi))$. The automorphism 
$\psi:\,\g\to\g$ induces an isometry from the orthogonal complement of 
$\h=\on{ker}(\Ad_{(\psi^{-1},g)}-1)$ onto $\mf{c}$.
Therefore, 
$${\det}_{\h^\perp}(\Ad_{(\psi^{-1},g)}-1)=
{\det}_{\mf{c}}(\Ad_g-\psi).$$
\end{proof}
\begin{proposition}\label{prop:b}
The Liouville measure of 
the symmetric space $G/G^\psi$ is equal to its Riemannian measure times 
a factor, $2^{\dim(G/G^\psi)/2}$.
\end{proposition}
\begin{proof}
In this case, $\mf{c}=\on{ker}(1-\psi)^\perp$ is the $-1$ 
eigenspace of $\psi$. Therefore 
$|{\det}_\mf{c}(1-\psi)|=2^{\dim(G/G^\psi)}$. 
\end{proof}
\begin{proposition}\label{prop:c}
The Liouville measure of $G$, viewed as a symmetric space for the
involution $\psi(a,b)=(b,a)$ of $G\times G$, coincides 
with its Riemannian measure. 
\end{proposition}
\begin{proof}
By Proposition \ref{prop:b}
the Liouville volume is equal to  $2^{\dim
G/2}\vol(G)^2$, divided by the Riemannian volume of the diagonal
subgroup $(G\times G)^\psi$. Since the induced Riemannian metric of 
$(G\times G)^\psi\cong G$ is twice the Riemannian metric of $G$, 
its volume is $2^{\dim G/2}\vol G$. 
\end{proof}

\begin{proposition}\label{prop:double}
The Liouville measure for the double $D(G)$ and for its fusion 
$\ti{D}(G)$ coincide with the Riemannian measure of $G\times G$. 
The Fourier coefficients of 
their DH-measures are 
\beq
\l \m_{D(G)},\chi_\lambda\otimes\chi_\mu\r
&=& \vol G^2 \ \delta_{\lambda,\mu}\\
\l \m_{\ti{D}(G)},\chi_\lam\r&=& \f{\vol G^2}{\dim V_\lam}\\
\eeq
\end{proposition}
\begin{proof}
The identification of Liouville measures and Riemannian measures 
follows from Proposition \ref{prop:c}, since fusion does not 
change the Liouville measure (Theorem \ref{th:fusion}).
Hence the DH-measure $\m_{D(G)}$ is the push-forward of the 
Riemannian measure on $G^2$ under
the moment map $(a,b)\mapsto(ab,a^{-1}b^{-1})$. Its Fourier
coefficients are obtained from the calculation, 
$$ 
\l \m_{D(G)},\chi_\lambda\otimes\chi_\mu\r
=\int_{G\times G} \chi_\lambda(ab)\chi_\mu(a^{-1}b^{-1})\d\vol_{G\times G}
=\vol G \ \delta_{\lambda,\mu}
$$
where we have used conjugation invariance and the orthogonality relations 
for irreducible characters. The formula for the DH-measure of 
$\ti{D}(G)$ follows from \eqref{eq:conv0}.
\end{proof}

\begin{proposition}
Let $\Co_1,\ldots,\Co_r$ be a collection of conjugacy
classes in $G$ and $g_j\in \Co_j$.  
The Liouville measure of the fusion product 
\begin{equation}\label{eq:modulispace}
 M=\ti{D}(G)\fus \cdots\fus \ti{D}(G)\fus\Co_1\fus\cdots\fus \Co_r
\end{equation}
(with $h\ge 0$ factors of $\ti{D}(G)$) is equal to the Riemannian
measure of $G^{2h}\times G/G_{g_1} \times\cdots\times G/G_{g_r}$, times a 
factor $\prod_{j=1}^r|{\det}_{\g_{g_j}^\perp}(\Ad_{g_j}-1)|^{1/2}$. 
The Fourier coefficients of its DH-measure $\m\in\ca{E}'(G)$ are given 
by the formula, 
\begin{equation}\label{eq:108}
 \l\m,\chi_\lam\r
=\vol G^{2h}\,\ 
\f{\prod_{j=1}^r \Vol(\Co_j)\,\chi_\lambda(\Co_j)}{
(\dim V_\lambda)^{2h+r-1}
}
\end{equation}
\begin{proof}
This follows directly from Propositions \ref{prop:LioConj}
and \ref{prop:double}, together with the formula \eqref{eq:conv} for the 
Fourier coefficients of convolutions of invariant measures.
\end{proof}
\end{proposition}

\section{Moduli spaces of flat bundles over surfaces}
\label{WittenFormulas}
\subsection{Connected semi-simple groups}\label{subsec:1}
Let $\Sig$ be a compact, oriented surface of
genus $h$ with $r$ boundary components, and let
$\underline{\Co}=(\Co_1,\ldots,\Co_r)$ be a collection of conjugacy
classes in $G$.  Denote by $\M_G(\Sig,\underline{\Co})$ the space of
isomorphism classes of flat $G$-bundles over $\Sig$, with holonomy
around the $j$th boundary circle in the conjugacy class $\Co_j$.  The
Atiyah-Bott \cite{at:ge} construction shows that
$\M_G(\Sig,\underline{\Co})$ is a compact symplectic space (sometimes
singular). $\M_G(\Sig,\underline{\Co})$ may also be obtained
\cite{al:mom} as a symplectic quotient of a group valued Hamiltonian
$G$-manifold
\begin{equation}
\label{eq:modreduction}
\M_G(\Sig,\underline{\Co})=M\qu G 
\end{equation}
where $(M,\om,\Phi)$ is the group valued Hamiltonian $G$-manifold defined 
in \eqref{eq:modulispace}. Recall that the $G$-action
on $M$ is given by conjugation on each factor in \eqref{eq:modulispace}.  
In particular, the
center $Z(G)\subset G$ acts trivially.  

In the following discussion we will make the assumption ($*$): Each
component of $\Phinv(e)$ contains at least one point with stabilizer
equal to $Z(G)$.  ($*$) implies that the number $k$ in
\eqref{eq:volform1} is equal to $\#Z(G)$.  ($*$) is automatic if $h
\ge 2$. Indeed $G$ can be generated\footnote{A closed subgroup
$H\subset G$ of a compact Lie group is {\em generated} by $S\subset G$
if it is the smallest closed subgroup containing $S$.}  by two
elements in $G$, and the commutator map $(a,b)\mapsto [a,b]$ is
surjective, see e.g. \cite[Corollary 6.56]{ho:st}, or Lemma
\ref{lem:surj} below).  This implies that the principal stabilizer for
the $G$-action on $G^2$ is equal to $Z(G)$ For surfaces of low genus
$h\le 1$, the validity of ($*$) depends on the conjugacy classes
$\Co_j$. Indeed, ($*$) never holds for $2h+r\le 2$. On the other hand,
one finds that $(*)$ holds for $h=1,r\ge 1$ if at least one of the
conjugacy classes $\Co_j$ is regular.  For a careful description of
the stabilizer groups in general, see Bismut-Labourie \cite[Section
5.5]{bi:sy}.

Using \eqref{eq:108} we obtain, 
\begin{equation}\label{eq:GeneralVolume}
\Vol(\M_G(\Sig,\ol{\Co}))=\#Z(G)\,\vol G^{2h-2}
\lim_{t\to 0+} \sum_{\lambda\in\on{Irr}(G)}e^{-t p(\lambda)}
\f{\prod_{j=1}^r \Vol(\Co_j) \chi_\lambda(\Co_j)}{(\dim
V_\lambda)^{2h+r-2}}
\end{equation}
which is Witten's formula \cite{wi:tw}. 
For $h\ge 2$, the sum is absolutely convergent 
at $t=0$, and the convergence factor may be omitted. 

\subsection{Connected components of the moduli space}
The connected components of the moduli space $\M_G(\Sig,\underline{\Co})$
are labeled by the topological type of the bundle. If $r\ge 1$, 
every $G$-bundle 
over $\Sig$ is isomorphic to the trivial bundle and therefore 
$\M_G(\Sig,\underline{\Co})$ is connected. 
Suppose $r=0$. Let $D\subset \Sig$ be a 
disk, and $\Sig'=\Sig\backslash D$. Any map $\gamma:\,\p D\to G$ defines 
a $G$-bundle over $\Sig$, by using $\gamma$ as a gluing function for 
trivial $G$-bundles over $D$ and $\Sig\backslash D$.  It is well-known 
that this construction gives a bijection between elements of $\pi_1(G)$ 
and isomorphism classes of $G$-bundles over $\Sig$. 

We compute the volumes of the connected components, for $h\ge 2$, as
follows.  Let $\pi:\,\tilde{G}\to G$ be the universal cover of $G$. For any 
$c\in \pi^{-1}(e)\subset Z(\ti{G})$ let 
$$\tilde{\Phi}^{(c)}: \,\ti{D}(G)\fus\ldots\fus\ti{D}(G) \to
\tilde{G}$$ 
the unique lift such that $\ti{\Phi}^{(c)}(e,\ldots,e)=c^{-1}$. The
quotient $\M_G^{(c)}(\Sig) =(\tilde{\Phi}^{(c)})^{-1}(e)/\tilde{G}$ is
the moduli space of flat connections on the bundle parametrized by
$c$.  We have
$$ \M_G(\Sig)=\coprod_{c\in \pi_1(G)} \M_G^{(c)}(\Sig).$$
The DH-measure $\dh$ with respect to  $\tilde{\Phi}^{(c)}$ 
has Fourier coefficients
$$ \l\dh^{(c)},\chi_\lambda\r=\vol G^{2h} \f{\chi_\lam(c^{-1})}{(\dim V_\lambda)
^{2h-1}},\ \ \lam\in \on{Irr}(\ti{G}).
$$
Since 
the principal stabilizer for the $\tilde{G}$-action has 
$\#Z(\tilde{G})=\#Z(G)\, \#\pi_1(G)$ elements, and  
$\Vol(\tilde{G})=\#\pi_1(G)\,\vol G$, Equation \eqref{eq:FourierDecomp} 
yields
\begin{equation}\label{eq:WittenFormula}
\Vol(\M_G^{(c)}(\Sig))=\vol G^{2h-2}
\f{\# Z(G)}{\#\pi_1(G)}
\sum_{\lambda\in\on{Irr}(\ti{G})}
\f{ \chi_\lambda(c^{-1})}
  { (\dim V_\lambda)^{2h-1}}.
\end{equation}
This agrees with Witten's result  (\cite{wi:tw}, Section 4.1).  
Summing over all $c\in \pi_1(G)$
we recover the formula for $\Vol(\M_G(\Sig))$ 
(as a sum over irreducible characters
for $G$ rather than $\ti{G}$).

\subsection{Two-component groups} 

If $G$ is disconnected, the orbit space $M/G$ may have several
connected components, with different principal stabilizer groups.  In
this section we consider moduli spaces of flat $G$-bundles over a
closed surfaces of genus $h \ge 2$ where $G$ is a compact Lie group
which has exactly two connected components and the identity component
$G^\o$ is semi-simple.  This class of groups is small enough so that
nice formulas exist, but large enough to include interesting examples, such
as moduli spaces of Riemannian (non-oriented) vector bundles on surfaces.
\begin{lemma} \label{lem:ztwo}
If $G$ has two components, the principal stabilizer for the $G$-action
on $G\times G$ is equal to the centralizer $Z_G(G^{\on{o}})$ on
$G^{\on{o}}\times G^{\on{o}}$, and equal to the center $Z(G)$ on the
other components of $G\times G$.
\end{lemma}

\begin{proof}
It suffices to show $G$ generated by elements $g_1\in G\backslash
G^{\on{o}}$ and $g_2\in G^{\on{o}}$.  Let $T$ be a maximal torus of
$G^{\on{o}}$, and $g_2\in T$ a regular element generating $T$.  Choose
$g_1\in G\backslash G^\o$ such that $\Ad_{g_1}(g_2)$ is not contained
in a proper closed subgroup $H\subset G^{\on{o}}$ containing $T$. This
is possible since $\{\Ad_{g_1}(g_2)|\, g_1\in G\backslash
G^{\on{o}}\}$ is a regular conjugacy class.  Then $g_2,\Ad_{g_1}(g_2)$
generate $G^{\on{o}}$ and therefore $g_1,g_2$ generate $G$.
\end{proof}

Let $G$ be a compact Lie group with two connected components. 
If $Z(G)\subseteq G^\o$, 
then $Z_G(G^\o)=Z(G^\o)$ while $Z(G)\subseteq 
Z(G^\o)$ is the subgroup fixed under the action of 
the component group $G/G^\o$. If $Z(G)\not\subseteq G^\o$, 
then $Z(G)=Z_G(G^\o)$ contains $Z(G^\o)$, and $Z(G)/Z(G^\o)=\Z_2$.

\begin{examples} 
\begin{enumerate}
\item
 \label{ex:A}
Let $G=\on{O}(n),\ n\ge 3$.  Then $Z(G)=Z_G(G^{\on{o}})=\{I,-I\}$. 
The spaces
$M_G(\Sig,\underline{\Co})$ are moduli spaces of 
flat Riemannian vector bundles.
\item
\label{ex:B}
Let $G=\on{O}(2)$. Then $Z(G)=\Z_2$ but $Z_G(G^{\on{o}})=\SO(2)$. 
\item
\label{ex:C}
Let $G^{\on{o}}$ be a compact, connected Lie group, and
$\psi\in\on{Aut}(G^{\on{o}})$ an involutive
automorphism. Let $G=\Z_2\ltimes G^{\on{o}}$ be the
corresponding semi-direct product. If $\psi$ is inner, 
$Z(G)=Z_G(G^{\on{o}})=\Z_2\times Z(G^{\on{o}})$. If $\psi$ is 
{\em not} inner, $Z_G(G^{\on{o}})=Z(G^{\on{o}})$ while 
$Z(G)=Z(G^{\on{o}})^\psi$. 
\item
\label{ex:D}
Let $G=\{A\in\U(n)|\,\det(A)^2=1\}$, $n\ge 2$. Note that $G$ is 
not a semi-direct product. 
One finds $Z(G)=Z_G(G^{\on{o}})=\{aI|\,a^{2n}=1\}=\Z_{2n}$.
The spaces $M_G(\Sig,\underline{\Co})$ are moduli 
spaces of flat Hermitian bundles over $\Sig$, such that the 
square of the canonical line bundle is trivial.
\end{enumerate}
\end{examples}

Let $\on{Irr}(G,G^{\on{o}})\subset \on{Irr}(G)$ be the set of
$G$-representations which are irreducible as
$G^{\on{o}}$-representations.  The image of the map
$\on{Irr}(G,G^{\on{o}})\to \on{Irr}(G^{\on{o}})$ consists of
irreducible representations that are fixed, up to isomorphism, by the
action of the component group $G/G^{\on{o}}$.
\begin{theorem}
Let $G$ be a compact Lie group with two components, with $G^{\on{o}}$
semi-simple, $\Sig$ be a closed surface of genus
$h\ge 2$,  $\M_G^{\on{o}}(\Sig) \subset \M_G(\Sig)$ be the equivalence
classes of flat  $G$-bundles that reduce to flat $G^{\on{o}}$-bundles, and 
$\M_G^1(\Sig)$ the remaining components. Then,
$$
\Vol(\M_G^{\on{o}}(\Sig))=\f{\#Z_G(G^{\on{o}})}{\#Z(G^{\on{o}})}
\Vol(\M_{G^{\on{o}}}(\Sig)),
$$
and 
$$  \Vol(\M_G^1(\Sig))=
(1-2^{-2h})\,\#Z(G)\,\vol G^{2h-2}
\sum_{\lambda\in \on{Irr}(G,G^{\on{o}})} (\dim
V_\lambda)^{2 -2h}. 
$$
\end{theorem}

Since we assume $h\ge 2$, the sum over $\lambda$ 
is absolutely convergent.

\begin{proof}
Let $M= \ti{D}(G)\fus \cdots\fus \ti{D}(G)$ 
be the fusion product of $h$ copies of $\ti{D}(G)$, and 
$\Phi:\,M\to G$ the moment map. By Lemma \ref{lem:lemma} from 
the Appendix B, all components of $\Phinv(e)$ meet 
the principal orbit type stratum of the corresponding component of $M$. 
Let  $M^o=\ti{D}(G^\o)\fus\cdots\fus \ti{D}(G^\o)$ be a fusion product
of $h$ copies of $\ti{D}(G^\o)$. The principal stabilizer is 
equal to $Z(G)$ on $M^1=M\backslash M^\o$ and to 
$Z_G(G^\o)$ on $M^\o$. An  open dense subset of
the symplectic quotient $M^{\on{o}}\qu G^{\on{o}}$
fibers over $M^{\on{o}}\qu G$ with fiber $Z_G(G^{\on{o}})/Z(G^{\on{o}})$.
Hence the volumes are related by a factor, 
$\f{\#Z_G(G^{\on{o}})}{\#Z(G^{\on{o}})}$. 

To compute $\Vol(M^1\qu G)$, we 
 calculate the Fourier coefficients $\l\m^1,\chi_\lam\r$ for $M^1$. 
Since $\chi_\lam$ vanishes on $G\backslash G^{\on{o}}$ for
$\lambda\not\in \on{Irr}(G,G^{\on{o}})$, we only need to consider 
$\lambda\in \on{Irr}(G,G^{\on{o}})$. We have 
$\l\m^1,\chi_\lam\r=\l\m,\chi_\lam\r-\l\m^{\on{o}},\chi_\lam\r$. 
The Fourier coefficients $\l\m,\chi_\lam\r$ are given by 
\eqref{eq:108} as usual. We claim that $\l\m^{\on{o}},\chi_\lam\r$ 
is given by a similar formula but with $G^\o$ in place of $G$: 
$$ \l\m^{\on{o}},\chi_\lam\r=2^{-2h} \vol G^{2h}\,\ 
(\dim V_\lambda)^{1-2h} $$
Indeed, $ \l\m_{\ti{D}(G^{\on{o}})},\chi_\lam\r=
\vol(G^\o)^2\,(\dim V_\lam)^{-1}$ by the formula for 
$G^\o$, since $\chi_\lam$ restricts to an irreducible character of 
$G^{\on{o}}$. It follows that 
$$\l\m^1,\chi_\lam\r=(1-2^{-2h})
\vol G^{2h}\,\ 
(\dim V_\lambda)^{1- 2h}
$$
These Fourier coefficients contribute with a factor $\#Z(G)$.
\end{proof}

%

\section{Mixed DH-distributions}
In this Section we generalize the definition of DH-measures for
Hamiltonian $G$-spaces $(M,\om,\Phi)$ with group valued moment map. We
will associate to $(M,\om,\Phi)$ invariant distributions on $G$ which
are smooth at regular values of $\Phi$ and which encode certain
intersection pairings on symplectic quotients. The results of this Section 
depend heavily on techniques developed in \cite{al:no}.  
In \cite[Theorem 5.2]{al:gr}, it is shown how to calculate 
the Fourier coefficients of the mixed DH-distributions by 
localization techniques.

\subsection{Mixed DH-distributions for Hamiltonian $G$-manifolds}
For Hamiltonian $G$-manifolds $(M,\om,\Phi)$ with $\g^*$-valued moment
maps, mixed DH-distributions were first introduced by Jeffrey-Kirwan
in \cite{je:lo1}, and studied in detail by Duistermaat \cite{du:eq},
Vergne \cite{ve:no} and Paradan \cite{pa:mo}.  The definition of mixed
DH-distributions uses the framework of equivariant cohomology. In
contrast to the above references, we will use the Weil model of
equivariant cohomology rather than the Cartan model.

As before we identify $\g^*\cong\g$ by a given invariant inner
product, and let $e_i$ be an orthonormal basis of $\g$. The Weil
algebra is the $W_G=S\g\otimes \wedge\g$, equipped with the $G$-action
induced by the adjoint action on $\g$. Let $L_i^{W}=L_i^{S\g}\otimes
1+1\otimes L_i^{\wedge\g}$ be the generators for the $G$-action,
$\iota_i^W=1\otimes \iota_i^{\wedge\g}$ the contraction operators, and
let
\begin{equation}\label{Differential}
\d^W=\sum_i (L_i^{S\g}\otimes 1)y_i+
\sum_i (v_i-\hh \sum_{j,k} f_{ijk} y_j y_k)\iota_i^W.
\end{equation}
Here $v_i,y_i$ are the generators of $S\g,\wedge\g$ corresponding to
the basis $e_i$, and $f_{ijk}=[e_i,e_j]\cdot e_k$ are the structure
constants.  Then $W_G$ together with these derivations is a
$G$-differential algebra. That is, the derivations $\iota_i^W,L_i^W$
and $\d^W$ satisfy (super-) bracket relations analogous to those for
contractions, Lie derivatives, and exterior differential for a
$G$-manifold $M$.  Given a $G$-manifold $M$, the equivariant
cohomology $H_G(M)$ is the cohomology of the basic subcomplex of
$W_G\otimes \Om(M)$, consisting of invariant elements that are
annihilated by contractions $\iota_i^W\otimes 1+1\otimes\iota_i$.  We
will also need the Weil algebra with distributional coefficients
$\wh{W}_G=\ca{E}'(\g)\otimes\wedge\g$, where $\ca{E}'(\g)$ is the
convolution algebra of compactly supported distributions on $\g$.  The
inclusion $S\g\hra \ca{E}(\g),\ \ v_i\mapsto \f{d}{d\,t}|_{t=0}
\delta_{t e_i}$ makes $W_G$ into a $G$-differential subalgebra of
$\wh{W}_G$.

Suppose now that $(M,\om_0,\Phi_0)$ is a Hamiltonian $G$-manifold. 
The {\em equivariant Liouville form} is a basic, closed element 
of $\wh{W}_G\otimes\Om(M)$ defined by 
\begin{equation}
\label{eq:L0} 
\L_0=e^{\om_0}\Phi_0^*\Lambda_0
\end{equation}
where 
\begin{equation}
\label{eq:Lambda0} 
\Lam_0=e^{-\sum_i y_i d\mu_i}\exp(-\hh \sum_{i,j,k} f_{ijk} y_i y_j \mu_k)\ \delta_\mu\in \wh{W}_G\otimes \Om(\g)
\end{equation}
Notice that the horizontal projection 
$P_{\on{hor}}^{\wedge\g}\L_0$ is just $e^{\om_0}\delta_{\Phi_0}$. 
Therefore, if $M$ is compact, the integral of 
$\L_0$ over $M$ coincides with  the DH-measure:
$$ \int_M \L_0
=\int_M P_{\on{hor}}^{\wedge\g}\L_0
=\int_M e^{\om_0}\delta_{\Phi_0}
=(\Phi_0)_*\nu_0
=\m_0.$$
Generalizing this equation, one defines for 
any equivariant cocycle ${\beta_0}\in (W_G\otimes \Om(M))_{basic}$, 
the {\em mixed DH distribution}
$\dh_0^{\beta_0}\in \ca{E}'(\g)^G$ by
\begin{equation}\label{eq:dh0}
 \dh_0^{\beta_0}=\int_M \beta_0\L_0.
\end{equation}
The map ${\beta_0}\mapsto \dh_0^{\beta_0}$ descends to a map 
$H_G(M)\to \ca{E}'(\g)^G$. We quote from \cite{je:lo1} the following 
properties of the mixed DH-distributions. $\dh_0^{\beta_0}$ is 
supported on the image of $\Phi_0$ and has 
singular support in the set of singular values of
$\Phi_0$. If $0$ is a regular value of $\Phi_0$, let the 
{\em Kirwan map}
$$\kappa:\,H_G(M)\to H_G(\Phi_0^{-1}(0))\cong H(M\qu G)$$
be the map given by restriction to the zero level set.
If $M/G$ is connected and each connected component of 
$\Phinv(0)$ meets the principal orbit type stratum, 
\begin{equation}\label{eq:hamilt}
\int_{M/\!/G} \kappa(\beta_0)\exp(\om_0)_{red}=
\f{k}{\vol G}\ \f{\dh^{\beta_0}_0}{\d\Vol_\g}\Big|_0.
\end{equation}
Here $(\om_0)_{red}$ the reduced symplectic form, and $k$ 
is the number of elements in a principal stabilizer.

\subsection{Mixed DH-distributions for group valued moment maps}
\label{NonComWeil}
To extend the construction of mixed DH-distributions to group-valued 
moment maps, we need a  non-commutative version of the 
Weil algebra introduced in \cite{al:no}. 
For simplicity, we will only discuss the case where $G$ is 
1-connected. The generalization to arbitrary compact 
groups (along the lines of Section \ref{sec:general}) is straightforward.

Let $U(\g)$ be the universal enveloping algebra, with generators 
$u_i$. It embeds into the convolution algebra $\ca{E}'(G)$
of distributions on $G$ by
the map $u_i\mapsto \f{d}{d\,t}|_{t=0} 
\delta_{\exp(t
e_i)}$. The non-commutative Weil algebra $\W_G$ and the non-commutative Weil
algebra $\wh{\W}_G$ with distributional coefficients are the 
$\Z_2$-graded algebras
$$\ca{W}_G:=U(\g)\otimes\on{Cl}(\g), \ \ 
\widehat{\ca{W}}_G:=\ca{E}'(G)\otimes \on{Cl}(\g).$$
The conjugation actions of $G$ on $\g$ and on $G$ induce an action on
$\W_G,\,\wh{\W}_G$. Let $L_i^\W$ be the corresponding Lie derivatives, and
$\iota_i^\W=\ad(x_i)$ (graded commutator with $x_i$).  As shown in
\cite{al:no},
$$ \d^\W=\ad(\sum_i u_i x_i -\f{1}{6}\sum_{i,j,k} f_{ijk} x_i x_j x_k)$$
is a differential, and makes $\W_G$ and $\wh{\W}_G$ into
$G$-differential algebras. For any $G$-manifold $M$, we define
$\wh{\ca{H}}_G(M)$ (resp.${\ca{H}}_G(M)$) 
as the cohomology of the basic subcomplex of
$\wh{\W}_G\otimes \Om(M)$ (resp. ${\W}_G\otimes \Om(M)$). Let 
$$ \La=e^{-\sum_i x_i\theta_i^R}\,\tau(g)\delta_g\in
\widehat{\W}_G\otimes \Om(G).$$
Now let $(M,\om,\Phi)$ be a group valued Hamiltonian $G$-manifold. From 
the properties of the form $\La$ (see \cite[Proposition 5.7]{al:no}) 
and the axioms for a space with group-valued moment map, it is immediate 
that the  {\em equivariant Liouville form}
$$ \L=e^\om\ \Phi^*\La \in \wh{\W}_G\otimes\Om(M)$$
is basic and closed. 
Note that if we denote
$$\Pi^\W=\int_G\otimes P_{\on{hor}}^\Cl:\,\wh{\W}^G
=\ca{E}'(G)\otimes \Cl(\g)\to \R$$ 
then  
$(\Pi^\W\otimes 1)\Lambda=\ca{T},\ \ (\Pi^\W\otimes 1)\L=\Gamma$, 
where 
$\ca{T}\in \Om(G)$ and $\Gamma\in\Om(M)$ are the differential forms 
introduced in Section \ref{subsec:cons}.
Again the integral of $\L$ over $M$ 
reproduces the DH-measure 
$\dh=\Phi_*(\Gamma_{[\dim M]})$:
$$ \int_M \L = 
\int_M P^{\on{Cl}}_{\on{hor}}(\L) 
=\int_M \Gamma \delta_\Phi=\Phi_*\nu=\dh.
$$
Given any equivariant cocycle 
$\beta\in (\W_G\otimes\Om(M))_{basic}$ we define 
\begin{equation}\label{eq:dhbeta}
\dh^\beta=\int_M \beta \L \in \ca{E}'(G)^G.
\end{equation}
The map $\beta\mapsto \dh^\beta$ vanishes on coboundaries, 
hence it descends to a map $\ca{H}_G(M)\to \ca{E}'(G)^G$.

%
%
%

\subsection{Exponentials}\label{Sec:Exponentials}
To study the properties of the mixed DH-distributions $\m^\beta$, we 
need to understand the relationship between the mixed DH-distributions 
of a Hamiltonian $G$-space $(M,\om_0,\Phi_0)$ and of its exponential. 
This involves the {\em quantization map} 
$\ca{Q}$ introduced in \cite{al:no}. Let 
$f(s)=\f{1}{s}-\hh \on{coth}(\f{s}{2})$, and let 
$T$ be the skew-symmetric tensor field on $\g\backslash
J^{-1}(0)$ given by $T(\mu)_{ij}=f(\ad_\mu)_{ij}$. 
Then 
\begin{equation}\label{eq:Q} 
\ca{Q}:\,\wh{W}_G\to \wh{\W}_G,\ \ \beta_0\mapsto 
\sig^{-1}\circ \exp_*\ 
\big(J^\hh \exp(\hh \sum_{i,j} T_{ij}\iota_i^W \iota_j^W)(\beta_0)\big)
\end{equation}
is a well-defined map, 
since the singularities of $T$ are compensated by the 
zeroes of $J^\hh$. It was shown in \cite{al:no}
that $\ca{Q}$ is a homomorphism of $G$-differential 
spaces, and for any $G$-manifold $M$ the induced map in 
cohomology $\wh{H}_G(M)\to \wh{\H}_G(M)$ is in fact a ring 
homomorphism. The quantization map $\ca{Q}$ has the property, 
\begin{equation}\label{Eq102}
\ca{Q}(\La_0)=e^\varpi (1\otimes \exp^*)\La.
\end{equation}
where $\varpi\in\Om^2(\g)$ is the 2-form from \ref{subsec:exp}.
Now let $(M,\om_0,\Phi_0)$ be a Hamiltonian $G$-space such 
that $\exp$ has maximal rank on $\Phi_0(M)\subset\g$, and 
let $(M,\om,\Phi)$ be its exponential. Let $\L_0,\L$ be 
the equivariant Liouville forms. 
Since $\om=\om_0+\Phi_0^*\varpi$, it follows immediately 
from \eqref{Eq102} that 
$$ \ca{Q}(\L_0)=\L.$$
Hence if $\beta_0\in (W_G\otimes\Om(M))$ is an equivariant cocycle 
and  $\beta=\ca{Q}(\beta_0)$, then $\ca{Q}(\beta_0\L_0)$ is 
cohomologous to $\beta\L$. It follows that 
\begin{equation}\label{Eq103}
 \m^\beta=J^\hh \exp_* (\m_0^{\beta_0}).
\end{equation}

\subsection{Intersection pairings on symplectic quotients}
Let $(M,\om,\Phi)$ be a compact group-valued Hamiltonian $G$-space. 
If $e$ is a regular value of $\Phi$ let 
$$ \kappa:\,\ca{H}_G(M)\to H(M\qu G)$$
be composition of the isomorphism $\ca{Q}^{-1}:\,\ca{H}_G(M)\cong H_G(M)$, 
pull-back to the level set $H_G(M)\to H_G(\Phinv(e))$, and 
isomorphism $H_G(\Phinv(e))\cong H(\Phinv(e)/G)$. 

\begin{theorem}
For every equivariant cocycle $\beta\in (\W_G\otimes\Om(M))_{basic}$  
the support of $\dh^\beta$ is contained in the image of $\Phi$, 
and the singular support in the set of singular values of $\Phi$.
If $e$ is a regular value of $\Phi$ and if 
each component of $\Phinv(e)$ meets $M_{\on{prin}}$, 
then 
\begin{equation}
\int_{M/\!/G} \kappa(\beta)
\exp(\om_{\on{red}})=
\f{k}{\vol G}\ \f{\dh^{\beta}}{\d\Vol_G}\Big|_e,
\end{equation}
using the notation of Theorem \ref{th:meaningofdh}. 
\end{theorem}
\begin{proof}
Write $P^{\on{Cl}}_{\on{hor}}(\beta)=\sum_J u^J\otimes \beta_J$ 
with $u^J\in U(\g)$ and $\beta_J\in \Om(M)$. 
Let $\odot$ be the following product structure on $\Om(M)$, 
$$ \gamma_1\odot\gamma_2=\diag_M^* \exp(-\hh \sum_i \iota_i^{1}
\iota_i^{2})(\gamma_1\otimes \gamma_2),$$
where $\diag_M:\,M\to M\times M$ is the diagonal embedding. 
Arguing as in the proof of Theorem \ref{th:fusion}, we have 
$P^{\on{Cl}}_{\on{hor}}(\beta\L) 
=P^{\on{Cl}}_{\on{hor}}(\beta)\odot P^{\on{Cl}}_{\on{hor}}(\L) 
=P^{\on{Cl}}_{\on{hor}}(\beta)\odot (\Gamma\delta_\Phi).$
Integrating over $M$, the terms involving contractions 
$\iota_i$ make no contribution, and we obtain 
$$
\int_M \beta\L=\int_M 
P^{\on{Cl}}_{\on{hor}} (\beta\L)=\int_M 
P^{\on{Cl}}_{\on{hor}}(\beta)\,\,\Gamma\delta_\Phi
 =\sum_J u^J \int_M 
(\beta_J \Gamma)\,\delta_\Phi
= \sum_J u^J \,\,\Phi_*(\beta_J\,\,\Gamma)_{[\dim M]}.
$$
Thus 
$$ \dh^\beta= 
\sum_J u^J \,\,\Phi_*(\beta_J\,\,\Gamma)_{[\dim M]}. 
$$
Since left convolution by $u^J$ is a differential operator, the
description of the support and singular support of $\dh^\beta$
follows. The interpretation of $\dh^\beta$ at the group unit $e$
follows as in the proof of Theorem
\ref{th:meaningofdh} by a reduction to the Hamiltonian case, using Equations
\eqref{eq:hamilt},\eqref{Eq103}.
\end{proof}

\begin{appendix}
\section{The spinning 4-sphere}
An interesting example of a Hamiltonian action with group valued moment map
is the 4-sphere $S^4$. We let $G=\SU(2)$, with inner product on 
$\mf{su}(2)$ given by $\xi_1\cdot\xi_2=-\hh
\on{tr}(\xi_1\xi_2)$. Cover $S^4$ by two coordinate charts
$U^\pm$ given as open balls $\Vert z \Vert^2<2$ in $\C^2$ 
with transition map
$$\varphi:\ U^+\backslash \{0 \} \to 
U^-\backslash \{ 0 \}, \ \ (z_1,z_2) \mapsto
(-\ol{z}_2,\ol{z}_1) \f{\sqrt{2-\Vert z\Vert ^2}}{\Vert z\Vert }.$$
The map $\varphi$ is equivariant for the standard
$\SU(2)$-action on $U^\pm\subset \C^2$. 
The $\SU(2)$-action on $S^4$ obtained in this way 
may be identified with the action induced from the embedding of 
$S^4$ in $\R^5 \cong \C^2 \oplus \R$, with $\SU(2)$ acting on the first factor
via the standard representation. Let $\Phi^\pm:\,U^\pm\to \SU(2)$
be given by
$$\Phi^+(z_1,z_2)=
\cos\big(\f{\pi \Vert z\Vert ^2}{2}\big)I + \f{i}{\Vert z\Vert ^2} \sin 
\big(\f{\pi \Vert z\Vert ^2}{2}\big)
\left(\begin{array}{cc} |z_1|^2-|z_2|^2 & 2 z_1\ol{z}_2\\ 
2 \ol{z}_1 z_2 & |z_2|^2 -|z_1|^2\end{array} 
\right)
$$
and $\Phi^-(z_1,z_2)=-\Phi^+(z_1,z_2)$. It is straightforward to
verify that $\varphi^*\Phi^-=\Phi^+$, so that $\Phi^\pm$ patch
together to define an $\SU(2)$-equivariant map $\Phi:\,S^4\to \SU(2)$.
Define  
$$ \om^+=\f{\pi}{2 \Vert z\Vert ^2}
\on{Re}(z_1\d \ol{z}_1+\d z_2\ol{z}_2)
\wedge \d\Vert z\Vert ^2
-\f{\sin(\pi \Vert z\Vert ^2)}{\Vert z\Vert ^4} \on{Im}
(\ol{z}_1 z_2 \d z_1 \d \ol{z}_2) 
$$
and let $\om^-=\om^+$. Again, one can check that 
$\om^\pm$ are $\SU(2)$-invariant and that $\varphi^*\om^-=\om^+$. 
\begin{theorem}
The triple $(S^4,\om,\Phi)$ is a group valued Hamiltonian $\SU(2)$-manifold. 
The volume form $\Gamma_{[4]}$ is given 
in terms of the standard volume form $\d\Vol_{\C^2}$ on $U^+\subset \C^2$ by 
$$\Gamma_{[4]} |_{U_+} =\f{\sin(\f{\pi}{2} \Vert z\Vert ^2)}{\f{\pi}{2}
\Vert z\Vert ^2} \d\Vol_{\C^2}.$$
\end{theorem}

\begin{proof}[Outline of proof]
We show that $(U^+,\om^+,\Phi^+)$ is the exponential (cf. Section
\ref{Sec:Exponentials}) of a Hamiltonian $\SU(2)$-manifold.  Let
$\om_0=\f{i}{2}(\d z_1 \d \ol{z}_1 +\d z_2 \d \ol{z}_2)$ be the
standard symplectic structure on $\C^2$. 
The defining $\SU(2)$-action on $\C^2$ is Hamiltonian, 
with moment map 
$$ \Phi_0(z_1,z_2)= \f{i\pi}{2} \left(\begin{array}{cc}
|z_1|^2-|z_2|^2 & 2 z_1\ol{z}_2\\ 2 \ol{z}_1 z_2 & |z_2|^2
-|z_1|^2\end{array} \right).
$$
To compute $\exp(\Phi_0)$ we use the formula 
$$ \exp(\xi)=\cos(\Vert \xi\Vert )I+\f{\sin(\Vert \xi\Vert )}{\Vert \xi\Vert } \xi, \ \
\xi\in\mf{su}(2)
$$
Putting $\xi=\Phi_0(z)$, and using that 
$\Vert \Phi_0(z)\Vert =\f{\pi \Vert z\Vert ^2}{2}$, 
we obtain $\exp(\Phi_0)=\Phi^+$. The symplectic 
form $\om_0$ can be written in the form 
$$ \om_0=\f{\pi}{2 \Vert z\Vert ^2}
\on{Re}(z_1\d \ol{z}_1+\d z_2\ol{z}_2)
\wedge \d\Vert z\Vert ^2
-\f{\pi}{\Vert z\Vert ^2} \on{Im}
(\ol{z}_1 z_2 \d z_1 \d \ol{z}_2) .
$$
%
%
The term $ -\f{1}{\Vert z\Vert ^4}\on{Im}(\ol{z}_1 z_2 \d z_1 \d \ol{z}_2)$ 
appearing in the formulas for both $\om^+$ and $\om_0$ 
is the pull-back of the normalized volume form $\d\Vol_{S^2}$ on 
$S^2$ under the quotient map $\psi:\,\C^2\backslash\{0\}\to \C P(1)$. 
Thus
$$\om^+-\om_0=(\sin(\pi \Vert z\Vert ^2)-\pi \Vert z\Vert ^2) \psi^*\d\Vol_{S^2}.$$
We want to identify this expression with $\Phi_0^*\varpi$. 
Introduce polar coordinates by the 
map $[0,\infty)\times S^2\to \mf{su}(2), (t,\zeta)\mapsto t\zeta $ 
where $S^2$ is viewed as the unit sphere in $\mf{su}(2)$. Then the 
form $\varpi$ is a multiple of the normalized volume form on $S^2$: 
$$ \varpi=\big(\sin(2 \Vert \xi\Vert )-2\Vert \xi\Vert \big)\ \d\Vol_{S^2}.$$
Using again $\Vert \Phi_0(z)\Vert =\f{\pi \Vert z\Vert ^2}{2}$, 
it follows that $\om^+=\om_0+\Phi_0^*\varpi$.
This shows that $(U^+,\om^+,\Phi^+)$ is the exponential of the subset
$\Vert z\Vert ^2<2$ of the Hamiltonian $\SU(2)$-manifold
$(\C^2,\om_0,\Phi_0)$. Multiplying the moment map $\Phi_+$ 
by the central element $-I\in \SU(2)$ one obtains 
$(U^-,\om^-,\Phi^-)$. 
The  volume form $\Gamma_{[4]}$ is obtained from 
Proposition \ref{prop:exp}, since $\Phi_0^* J^{1/2}
=\f{\sin(\Vert \Phi_0\Vert )}{\Vert \Phi_0\Vert }$. 
\end{proof}
The spinning 4-sphere is an example of a multiplicity free space: 
All reduced spaces for the $\SU(2)$-action are points. 
The $\SU(2)$-action extends to an action of $\U(2)$, by letting  
the central $\U(1)$ act on $U^+\subset \C^2$ by
$e^{i\phi}\cdot(z_1,z_2)=(e^{i\phi}z_1,e^{i\phi}z_2)$ 
and on $U^-$ by the opposite action, 
$e^{i\phi}\cdot(z_1,z_2)=(e^{-i\phi}z_1,e^{-i\phi}z_2)$. 
A moment map for the $\U(2)$ action is given by 
$$\ti{\Phi}^+(z_1,z_2)= \ti{\Phi}^-(z_1,z_2)= 
e^{-\f{i \pi}{2} \Vert z \Vert^2}\Phi^+(z_1,z_2).$$
The reduced spaces by the $\U(1)$-action are conjugacy classes of 
$\SU(2)$. In Hurtubise-Jeffrey \cite{hu:re}, the $S^4$ example appears 
as the {\em imploded cross-section} for the double $\ti{D}(\SU(2))$.

\section{Principal Stabilizers}\label{app:b}
The following Lemma is an extension of a well-known fact 
about compact connected semi-simple Lie groups
(Goto's commutator theorem, \cite[Theorem 6.55]{ho:st}):

\begin{lemma}\label{lem:surj}
Suppose $G$ is a compact Lie group, with $G^{\on{o}}$ semi-simple, and
assume that $G/G^\o$ is finite cyclic. Then the restriction of
the commutator map $f:\,G\times G\to G,\ (a,b)\mapsto [a,b]$ to any
component of $G\times G$ is onto $G^\o$.
\end{lemma}

\begin{proof} 
Note first that since the component group $G/G^\o$ is abelian, the
commutator map takes values in $G^\o$. Let $xG^\o\in G/G^\o$ be a
generator.  We want to show that $f$ is onto $G^\o$ on the component
of $(x^{k},x^{l})$, for any $k,l\in \Z_{\ge 0}$.  Using
$[a,b]=[b,a]^{-1}$, we may assume $k\ge l$. If $l>0$, let $m$ be the
largest non-negative integer such that $k-ml\ge 0$.  Using
$[a,b]=[ab^{-m},b]$ we see that $f$ is onto $G^\o$ on the component of
$(x^{k},x^{l})$ if and only it is onto $G^\o$ on the component of
$(x^{k-ml},x^{l})$.  Let $k_1=l,l_1=k-ml$. Iterating this procedure,
we obtain a finite sequence $k_j\ge l_j\ge 0$ such that $f$ is onto
$G^\o$ on the component of $(x^{k_j},x^{l_j})$. The sequence
terminates when $l_j=0$. We have thus reduced the problem to the case
$l=0$. Since $f$ is equivariant, it suffices to show that some maximal
torus $T\subset G^\o$ is in the image for this component.  Lemma
\ref{lem:dyn} below shows that one can find $a\in N_{G^\o}(T) x^k$
such that $\Ad_a$ preserves $\t$ and ${\Ad_a}|_\t$ has no eigenvalue
equal to $1$. Given $t\in T$, choose $\xi\in\t$ with
$\exp((\Ad_a-1)\xi)=t$. Then $b=\exp\xi$ satisfies
$$ [a,b]=\exp(\Ad_a\xi)\exp(-\xi)=
\exp(\Ad_a\xi-\xi)=t.$$
\end{proof}

\begin{lemma}\label{lem:dyn}
Let $\g$ be a compact, semi-simple Lie algebra, 
and $\t\subset\g$ a maximal abelian subalgebra. For any 
Lie algebra automorphism $\psi\in \on{Aut}(\g)$, 
there exists an inner automorphism $\phi\in \on{Int}(\g)$ 
such that $\phi\circ\psi$ preserves $\t$, and 
its restriction to $\t$ has no eigenvalue equal to $1$. 
\end{lemma}

\begin{proof}
Let $\t_+\subset \t$ be a positive Weyl chamber. By composing 
$\psi$ with an inner automorphism, we may assume 
that $\psi$ takes $\t_+$ to itself, and in fact any element 
of $\on{Out}(\g)=\on{Aut}(\g)/\on{Int}(\g)$ has a unique 
representative of this form. (That is, $\psi$ is induced
from an automorphism of the Dynkin diagram of $\g$.) We have to find 
a Weyl group element $w\in W$ such that $w\circ \psi|_\t$ has 
no eigenvalue equal to 1. Then any $\phi\on{int}(\g)$ preserving 
$\t$ and acting as $w$ on $\t$ will have the required property. 
If $\psi=1$, one can take $w$ to be any Coxeter element, 
cf. \cite[Chapter 3.16]{hu:re}. 
Suppose $\psi \not=1$. By decomposing $\g$ with respect to
the action of $\psi$, we may assume that $\g$ contains no proper
$\psi$-invariant ideal. We consider two 
cases: 

1) Assume that  $\g$ is not simple. Let 
$\h$ be a simple ideal of $\g$. Then $\g=\oplus_{j=0}^{k-1} 
\psi_0^j(\h)$, where $k$ is the order of $\psi_0$. Let $w'$ be a
Coxeter element of $\h$. Then $w=(w',1,\ldots,1)$ has the required 
property. 

2) We are left with the case of  $\g$ simple. 
Recall that $\on{Out}(\g)=\Z_2$ for the Lie algebras $A_n\ (n\ge
2),\, D_n\ (n\ge 5)$ and $E_6$, $\on{Out}(\g)=S_3$ for $D_4$, and
$\on{Out}(\g)=\{1\}$ in all other cases.

For the Lie algebras $A_n$ and $E_6$, and $\psi\not=1$, 
one can take $w$ to be the 
longest element of the Weyl group. Indeed $w\circ \psi|_\t=-\on{id}_\t$ 
in both cases (see \cite{bo:li}, planche I (p.251), 
planche V (p.261)). For the Lie algebra
$D_n$ ($n\ge 5$), consider the standard presentation of the root system as the 
set of lattice vectors in $\R^n$ of length $\sqrt{2}$. In the 
standard basis $\eps_i$ of $\R^n$, a system of simple roots 
is given by $\alpha_1=\eps_1-\eps_2,\alpha_2=\eps_2-\eps_3,\ldots,
\alpha_{n-1}=\eps_{n-1}-\eps_n,\,\alpha_n=\eps_{n-1}+\eps_n$.
The automorphism $\psi$ exchanges $\alpha_{n-1}$ and $\alpha_n$, 
hence $\psi|_\t$ is the linear map changing the sign of the last 
coordinate in $\R^n$. Take $w\in W$ to be a cyclic permutation of the 
coordinates. Then $w\circ \psi$ has no fixed vector. It remains to consider
the case of $D_4$. If $\psi$ is induced by a diagram automorphism of 
order $3$, the map $\psi|_\t$ has eigenvalues the third roots of unity.
The longest Weyl group element $w$ acts as $-\on{id}$ on $\t$, 
hence $w\circ \psi|_\t$ has no eigenvalue equal to $1$. Finally, if 
$\psi$ is induced by the automorphism exchanging $\alpha_3,\alpha_4$, 
we can take $w$ to be cyclic permutation of the coordinates as above, 
the other diagram automorphisms of order $2$ are obtained from this 
by conjugation with a third order automorphism. 
\end{proof}

\begin{lemma}\label{lem:lemma}
Let $G$ be a compact Lie group with two components and
with $G^\o$ semi-simple. 
For all $h\ge 2$, every connected component of every fiber 
of the map $\Phi:\,G^{2h}\to G,\,(a_1,b_1,\ldots,a_h,b_h)
\mapsto \prod_{j=1}^h [a_j,b_j]$ meets the principal 
orbit type stratum $(G^{2h})_{\on{prin}}$ for the conjugation action of $G$.
\end{lemma}

\begin{proof}
For $j=1,\ldots,h$ let $p_j:\,(G^2)^h\to G^2$ denote projection to the
$j$th factor. Given a component $X\subset G^{2h}$ it is possible to
choose $j$ such that the principal stabilizer for $p_j(X)$ equals that
of $X$. (The index $j$ is arbitrary if $X=(G^\o)^{2h}$; otherwise pick
$j$ such that $p_j(X)\not= (G^\o)^2$.) Thus $x\in X_{\on{prin}}$
whenever $p_j(x)\in (p_j(X))_{\on{prin}}$. It hence suffices to show
that for any given $y\in G^2$, $g\in G$, the fiber $p_j^{-1}(y)$ meets
each component of $\Phinv(g)$.  Interpret $G^{2h}$ as a fusion
product of $h$ copies of the double $\ti{D}(G)$, with $\Phi$ as its
moment map. By Lemma \ref{lem:surj}, the restriction of $\Phi$ to any
component $X\subset G^{2h}$ is surjective, and by another application
of this Lemma each fiber of $\Phinv(g)\cap X$ meets $p_j^{-1}(y)$. If
$G^\o$ is 1-connected, this completes the proof since each $X$
contains a {\em unique} component of $\Phinv(g)$. (Recall
\cite[Theorem 7.2]{al:mom} that the fibers of the moment map for a
compact, connected, group-valued Hamiltonian space are connected,
provided the group is 1-connected.)  If $G^\o$ is not 1-connected, we construct a finite cover $\hat{G}\to G$, with identity
component $\wh{G}^\o$ the universal cover of $G^\o$: Choose $g_0\in G$
such that $\psi:=\Ad_{g_0}:\,\g\to \g$ preserves a maximal abelian
subalgebra $\t\subset \g$. Any such $\psi$ has finite order $k>0$,
hence it defines an action of $\Z_k$ on $\hat{G}^\o$ by automorphisms.
We define $\hat{G}=\Z_k\ltimes
\hat{G}^\o$, with covering map $\hat{G}\to G,\
(\psi^{-l},\hat{h})\mapsto g_0^{-l}h$ where $h\in G^\o$ is the image of 
$\hat{h}\in\hat{G}^\o$. Let $\pi:\,\hat{G}^{2h}\to G^{2h}$ be the projection.  The product of commutators,
$\hat{\Phi}:\,\hat{G}^{2h}\to \hat{G}^\o$ is the moment map for the
$\hat{G}^\o$-action. The fiber $\Phinv(g)$ is covered by the union of
$\hat{\Phi}^{-1}(\hat{g})$, as $\hat{g}\in\hat{G}$ ranges over
pre-images of $g$. Again, Lemma \ref{lem:surj} shows that the
restriction of $\hat{\Phi}$ to any component $\hat{X}\subset
\hat{G}^{2h}$ is surjective, and $\hat{\Phi}^{-1}(\hat{g})\cap
\hat{X}$ meets $(p_j\circ \pi)^{-1}(y)$.  This proves the Lemma since
$\hat{\Phi}^{-1}(\hat{g})\cap \hat{X}$ is connected.
\end{proof}

\end{appendix}

\bibliographystyle{amsplain}

\providecommand{\bysame}{\leavevmode\hbox to3em{\hrulefill}\thinspace}
\providecommand{\MR}{\relax\ifhmode\unskip\space\fi MR }
\providecommand{\MRhref}[2]{%
  \href{http://www.ams.org/mathscinet-getitem?mr=#1}{#2}
}
\providecommand{\href}[2]{#2}

\end{document}